\documentclass[11pt]{article}
\usepackage[a4paper,margin=25mm]{geometry}
\usepackage{amsmath,amssymb,amsthm,mathtools,booktabs,hyperref}
\hypersetup{colorlinks=true,linkcolor=blue,citecolor=blue,urlcolor=blue}
\numberwithin{equation}{section}
\newtheorem{theorem}{Theorem}[section]
\newtheorem{lemma}[theorem]{Lemma}
\newtheorem{proposition}[theorem]{Proposition}

\theoremstyle{remark}\newtheorem{remark}[theorem]{Remark}
\newcommand{\R}{\mathbb R}
\newcommand{\dd}{\,\mathrm d}
\newcommand{\curl}{\operatorname{curl}}
\newcommand{\diver}{\operatorname{div}}

\newcommand{\Aminus}{\mathcal A_{-1}}
\title{Symmetry and uniqueness of degree~$-1$ minimizers \\
of a chiral magnetic energy}
\author{Xiuqing Duan\\
  {\small Division of Mathematical Sciences, Nanyang Technological University, Singapore}\\
  {\small\texttt{xiuqing.duan@ntu.edu.sg}}}
\date{}
\begin{document}
\maketitle
\begin{abstract}
We consider the two-dimensional energy of a chiral magnet, 
consisting of the exchange energy, a bulk Dzyaloshinskii--Moriya interaction with constant $\kappa\ne0$, 
and the Zeeman energy in an applied field $h>0$. 
We prove that the minimizers of this energy among maps of degree~$-1$ are exactly the translates of one axially symmetric Bloch skyrmion, 
whose helicity is determined by the sign of $\kappa$. 
This holds for $h \geq c \kappa^2$, where $c \in(\pi^2 / 16,1)$ is determined by the slope of the radial profile at the origin; 
in particular, it holds in the whole range $h\ge\kappa^2$, 
which answers a question raised by Melcher and by D\"oring and Melcher.
The proof is a calibration argument. 
Using the profile of the skyrmion, we write the energy density of an arbitrary map as a weighted topological density, 
plus a divergence, plus a nonnegative remainder, 
and we show that this identity can be integrated for every map of finite energy. 
The method also gives the uniqueness of the monotone solution of the profile equation for $h\ge\kappa^2/2$.
\end{abstract}
\tableofcontents

\section{Introduction}

Magnetic skyrmions are vortex-like patterns in the magnetization of a ferromagnet that carry a nonzero topological degree. 
In materials without inversion symmetry they are stabilized by the Dzyaloshinskii--Moriya interaction (DMI). 
They were predicted by Bogdanov and Yablonskii~\cite{BY89}, and isolated skyrmions were analysed in detail by Bogdanov and Hubert~\cite{BH94}. 
Twenty years later, a lattice of skyrmions was found by neutron scattering in the chiral magnet MnSi~\cite{Muh09}, 
skyrmions were imaged in real space in a thin film~\cite{Yu10}, and single skyrmions were written and deleted one at a time in an ultrathin film~\cite{Rom13}. 
Their small size and their stability have made them candidates for carriers of information in magnetic devices~\cite{FCS13}; see~\cite{NT13,FRC17} for reviews.

The simplest continuum model describes the magnetization of a thin film by a map $m\colon\R^2\to\mathbb S^2$ with energy
\begin{equation}\label{originalenergy}
E_{\kappa,h}(m)=\frac12\int|\nabla m|^2+\kappa D(m)+\frac h2\int|m-e_3|^2,
\end{equation}
the sum of the exchange energy, the bulk DMI with constant $\kappa\ne0$, 
and the Zeeman energy in an applied field $he_3$, $h>0$; 
note that $\frac12|m-e_3|^2=1-m_3$. 
Here and below, integrals without a domain are over $\mathbb R^2$. 
Formally $D(m)=\int m\cdot\operatorname{curl}m$, 
where, for fields independent of $x_3$,
\[
\operatorname{curl}m=(\partial_2m_3,\,-\partial_1m_3,\,\partial_1m_2-\partial_2m_1).
\]
We minimize $E_{\kappa,h}$ in the class
\[
\mathcal A_{-1}=\big\{m\in H^1_{\rm loc}(\mathbb R^2;\mathbb R^3):\ |m|=1\text{ a.e.},\
m-e_3\in H^1(\mathbb R^2;\mathbb R^3),\ Q(m)=-1\big\},
\]
where
\[
Q(m)=\frac1{4\pi}\int m\cdot(\partial_1m\times\partial_2m)\,dx
\]
is the degree; its integrand is integrable since $|m\cdot(\partial_1m\times\partial_2m)|\le\frac12|\nabla m|^2$.

For a general $m$ with $m-e_3\in H^1$, the density $m\cdot\operatorname{curl}m$ need not be integrable, 
so we define $D$ by continuity. 
Write $m=(w,m_3)$ with $w=(m_1,m_2)\in\mathbb R^2$, and let $J(x_1,x_2)=(-x_2,x_1)$. 
If $m-e_3\in C^\infty_c(\mathbb R^2;\mathbb R^3)$, then $\int e_3\cdot\operatorname{curl}(m-e_3)=0$, 
and integration by parts gives
\begin{equation}\label{dmiext}
\int m\cdot\operatorname{curl}m=\int(m-e_3)\cdot\operatorname{curl}(m-e_3)
=-2\int w\cdot J\nabla m_3 .
\end{equation}
The last two integrals are continuous with respect to $m-e_3$ in $H^1$. 
Hence, they agree whenever $m-e_3\in H^1$, and we define $D(m)$ as their common value.
Since $\int|\nabla\xi|^2-\int|\curl\xi|^2=\int(\partial_1\xi_1+\partial_2\xi_2)^2\ge0$ for $\xi\in H^1(\R^2;\R^3)$, 
\eqref{dmiext} gives $|D(m)|\le\|m-e_3\|_{L^2}\|\nabla m\|_{L^2}$. 
Consequently, as observed by Melcher~\cite{Mel14}, $E_{\kappa,h}\ge0$ if $h\ge\kappa^2$, 
and if $h>\kappa^2$ the energy is bounded from below by a positive multiple of
$\|\nabla m\|_{L^2}^2+\|m-e_3\|_{L^2}^2$. 
Theorem~\ref{originaltheorem} covers the range $h \geq \kappa^2$, including the endpoint $h=\kappa^2$, where the energy is nonnegative but the quadratic lower bound degenerates.

The mathematical analysis of \eqref{originalenergy} was started by Melcher~\cite{Mel14}. 
He proved by concentration-compactness that, for every $\kappa\ne0$ and $h>\kappa^2$, 
the energy attains its infimum on $\Aminus$, without any symmetry assumption. 
He also pointed out that it is a delicate question whether a minimizer among rotationally symmetric fields also minimizes the energy among fields of arbitrary shape. 
D\"oring and Melcher~\cite{DM17} extended the existence result, for small DMI, 
to the potentials $\int|m-e_3|^p$ with $2\le p\le4$, and studied the conformal limit, 
in which the DMI and the potential carry a small weight; 
for $2<p<4$ they showed that almost minimizers are then close to the axially symmetric Belavin--Polyakov profile~\cite{BP75}. 
For $p=4$ the energy has a Bogomol'nyi structure~\cite{Bog76}, 
and an explicit axially symmetric minimizer is available~\cite{DM17}. 
D\"oring and Melcher asked whether the minimizers in $\Aminus$ are unique up to translation and axially symmetric. 
The question has come up repeatedly since then. 
Li and Melcher~\cite{LM18} proved that for large $h$ the axisymmetric skyrmion is a strict local minimizer up to translation, 
and described its global minimality as a major open question. 
Gustafson and Wang~\cite{GW21}, who constructed co-rotational skyrmions close to harmonic maps, 
noted that it was not known whether the global minimizer has this symmetry.

For the quartic potential, the question has recently been answered: Ibrahim, Miura, Rom\'an and Shimizu~\cite{IMRS25} 
proved that below the critical coupling, the explicit Bogomol'nyi skyrmion is the only minimizer up to translation. 
At the critical coupling, the Bogomol'nyi equations have many solutions, 
some of which are not axially symmetric~\cite{BRS20,IMRS25}, and beyond it, the skyrmion is unstable, 
and the energy is unbounded from below~\cite{IS23}. 
For the Zeeman potential, which is the case of main physical interest, much less is known. 
Ibrahim and Shimizu~\cite{IS} constructed equivariant skyrmions for all values of the parameters by a global, 
non-perturbative argument, and proved their exponential decay, their monotonicity for small coupling, and their linear stability for small DMI.
For the related easy-axis anisotropy model, asymptotic descriptions of the profiles of small and large skyrmions were given in~\cite{KMV20,KMV21}. 
For ultrathin films, where the DMI is interfacial and the stray field is taken into account, 
Bernand-Mantel, Muratov, and Simon~\cite{BMS21} gave a quantitative description of the minimizers when the skyrmions are small, 
together with a rigidity estimate for degree~$\pm1$ harmonic maps; 
see also~\cite{MMSS23,MSS25} for skyrmions under confinement and minimizers of higher degree. 
None of these results gives the symmetry of minimizers for the Zeeman potential.

In this paper, we answer these questions for the Zeeman potential in the whole range $h>\kappa^2$ of Melcher's existence theorem and at its endpoint $h=\kappa^2$.

\begin{theorem}[Symmetry and uniqueness]\label{originaltheorem}
For every $\kappa\ne0$ and $h\geq \kappa^2$, put
$\nu=|\kappa|/\sqrt h\in(0,1]$. There is a unique regular monotone
profile $f_\nu$ with the endpoints and equation in \eqref{profileode}.
Every global minimizer of \eqref{originalenergy} on $\Aminus$ agrees
almost everywhere with exactly one field in the family
\begin{align}
 m_{\kappa,h,a}(x)
  ={}&\operatorname{sign}(\kappa)
       \sin f_\nu(\sqrt h\,|x-a|)
                   \left(\frac{J(x-a)}{|x-a|},0\right) \notag\\
   &+\cos f_\nu(\sqrt h\,|x-a|)e_3,
                           \qquad a\in\R^2,       \label{originalfamily}
\end{align}
with the value $-e_3$ at $x=a$, and every $m_{\kappa,h,a}$ is a minimizer. Thus
\[
 \arg\min_{\Aminus}E_{\kappa,h}
       =\{m_{\kappa,h,a}:a\in\R^2\}.
\]
Each $m_{\kappa,h,a}$ is smooth, takes the value $-e_3$ only at $x=a$, and satisfies
\[
 m_{\kappa,h,a}(a+R_\theta y)=\widehat R_\theta\, m_{\kappa,h,a}(a+y) \qquad
(\theta\in\mathbb R,\ y\in\mathbb R^2),
\]
where $R_\theta$ is the rotation of $\mathbb R^2$ by the angle $\theta$, 
and $\widehat R_\theta$ is the rotation of $\mathbb R^3$ by $\theta$ about $e_3$.
\end{theorem}

\begin{theorem}[The range of the calibration]\label{rangetheorem}
There is a number $\mu_{\rm cal}\in(1,4/\pi)$ with the following properties. 
For every $0<\mu\le\sqrt2$, the profile equation \eqref{profileode} with parameter $\mu$ has exactly one solution $f_\mu$;
its slope $a(\mu)=-f_\mu'(0)$ is a continuous, strictly decreasing function of $\mu$, and $a(\mu)\ge1$ if and only if $\mu\le\mu_{\rm cal}$. 
For every $\kappa\ne0$ and $h>0$ with $|\kappa|/\sqrt h\le\mu_{\rm cal}$, 
that is, $h\ge\mu_{\rm cal}^{-2}\kappa^2$, the conclusions of Theorem~\ref{originaltheorem} hold. 
On the other hand, if $h<\pi^2\kappa^2/16$, then $E_{\kappa,h}$ is unbounded from below on $\Aminus$.
\end{theorem}

No symmetry is assumed for the competitors in Theorem~\ref{originaltheorem}; 
the axial symmetry of the minimizers is part of the conclusion. 
Thus, the minimizer is unique up to translation, and its helicity is determined by the sign of $\kappa$. 
In particular, the minimizers found by Melcher~\cite{Mel14} are translates of a Bloch skyrmion, 
and the axisymmetric skyrmion is a global minimizer for every $h\ge \kappa^2$, 
which answers the question of Li and Melcher~\cite{LM18}. 
Since the lower bound is attained by a Bloch field, our argument also gives the existence of minimizers without concentration-compactness. 
In particular, minimizers exist for $h=\kappa^2$, which is not covered by~\cite{Mel14}.
To our knowledge, Theorem~\ref{originaltheorem} is the first result on the symmetry of minimizers for the Zeeman potential.

The uniqueness of $f_\nu$ is not assumed in the proof; it follows from the classification of minimizers. 
For large $h$, it was proved by Li and Melcher~\cite{LM18} through an analysis of the profile equation, 
and Ibrahim and Shimizu conjectured that the solutions of the profile equation lying below the Belavin--Polyakov profile are unique~\cite{IS}. 
Theorem~\ref{originaltheorem} gives uniqueness among monotone profiles for every $h\geq \kappa^2$, 
and Theorem~\ref{rangetheorem} for every $h\ge\kappa^2/2$.

We add two remarks. 
First, in two dimensions a rotation of the target by $\pi/2$ about $e_3$ 
changes the bulk DMI into the interfacial DMI $\int(m_3\diver w-w\cdot\nabla m_3)$, up to the sign of $\kappa$, 
and leaves the other terms in \eqref{originalenergy} and the degree unchanged. 
Hence Theorem~\ref{originaltheorem} also classifies the minimizers for the interfacial DMI, 
which is the relevant interaction in ultrathin films~\cite{Rom13,FRC17}: they are the translates of a N\'eel skyrmion. 
Second, the condition $h\ge\kappa^2$ enters the proof only through the bound $a_0=-f'(0)>1$ 
of Lemma~\ref{originslope}: the profile exists for $h\ge\kappa^2/2$ (Appendix~\ref{ex-existence}), 
and Sections~\ref{sec:target}--\ref{sec:equality} apply to every $\kappa\ne0$, $h>0$, and every profile with $a_0\ge1$ (Theorem~\ref{conditional}). 
Theorem~\ref{rangetheorem} makes the range of the method precise. 
Numerically, $a_0\ge1$ holds exactly for $h/\kappa^2\ge0.729$, and the energy of the Bloch skyrmion is negative for $h/\kappa^2<0.80$, 
the field below which the skyrmion lattice has lower energy than the uniform state~\cite{BH94}; these values are not used in the proofs.

\paragraph{Strategy and new ingredients.}
Standard tools for proving symmetry, such as rearrangement or the moving plane method, 
do not seem to apply here: the problem is vector valued, the constraint $Q(m)=-1$ is topological, 
and the DMI distinguishes a handedness. Our proof is a calibration argument instead. 
Its model is the Bogomol'nyi argument~\cite{BP75,Bog76} for the exchange energy alone:
\[
\frac12|\nabla m|^2=-m\cdot(\partial_1m\times\partial_2m)+\frac12|\partial_1m-m\times\partial_2m|^2,
\]
so the exchange energy of a map of degree $-1$ is at least $4\pi$, with equality exactly when $\partial_1m=m\times\partial_2m$. 
Melcher~\cite{Mel14} combined the exchange and DMI terms into helical derivatives 
$\partial_\alpha m-\kappa\,e_\alpha\times m$ 
and used a completion of the square of this type to bound the energy from below in sectors of nonnegative degree, 
and for the quartic potential a Bogomol'nyi bound is attained in degree $-1$~\cite{DM17}. 
For the Zeeman potential in degree $-1$, however, no such identity attained by the skyrmion was known.

By scaling and rotating the target, it suffices to consider $h=1$ and $\kappa=\mu>0$, and we write $E_\mu=E_{\mu,1}$. 
Our main observation is that the profile $f$ of the expected minimizer can be used to build an identity of Bogomol'nyi type. 
From $f$ we construct two functions $H$ and $F$ on $[0,\pi]$ and show that, 
writing $m=(\sin t\,n,\cos t)$ with $|n|=1$, the energy density of an arbitrary map equals
\[
H(t)\,m\cdot(\partial_1m\times\partial_2m)+\diver\big(F(t)\,Jn\big)+\mathcal R
\]
at every point where $m\ne\pm e_3$ (Lemma~\ref{lem:calibration}). 
The remainder $\mathcal R$ is a quadratic form in the difference between $\nabla m$ 
and the gradient of the Bloch field at the same polar angle, and it is positive definite because $|H|<1$. 
This strict bound is the heart of the matter. It is equivalent to the fact that $-f'(r)$, $\sin f(r)/r$, 
and $2\sin(f(r)/2)$ are the side lengths of a nondegenerate triangle for every $r>0$ (Lemma~\ref{triangles}); 
as $r\to0$, it reduces to $a_0\ge1$, and this is where the condition $h\ge\kappa^2$ enters the calibration, through Lemma~\ref{originslope}.

The second ingredient is the integrated form of this identity. For every $m\in\Aminus$ we show that
\[
E_\mu(m)=E_\mu(B_f)+\int_\Omega\mathcal R\dd x+2\,|\{m=-e_3\}|,\qquad
\Omega=\{m\ne\pm e_3\},
\]
where $B_f$ is the Bloch field with profile $f$ (Lemma~\ref{globalidentity}). 
Two points require care. 
The integral of the first term is given by the weighted degree formula 
$\int G(t)\,m\cdot(\partial_1m\times\partial_2m)=2\pi Q(m)\int_0^\pi G(t)\sin t\dd t$, 
which we prove for every continuous $G$ and every map of finite energy (Lemma~\ref{weighteddegree}), 
and the divergence term integrates to zero (Lemma~\ref{globalidentity}), although $m$ need not decay pointwise. 
The identity is exact. It gives the sharp lower bound $E_\mu(m)\ge E_\mu(B_f)$, 
and it shows that the energy excess of $m$ controls the $L^2$ distance between $\nabla m$ 
and the gradient of the Bloch field at the same polar angle; see \eqref{comparison}. 
In the equality case, this distance vanishes, 
so $m$ solves a first-order system whose solutions of degree $-1$ are exactly the translates of $B_f$ (Proposition~\ref{equality}).

The main new contributions of this paper are the pointwise identity of Lemma~\ref{lem:calibration}, 
the triangle inequalities of Lemma~\ref{triangles} that make it coercive whenever $a_0\ge1$, and its exact integrated form in Lemma~\ref{globalidentity}.
It would be interesting to know whether calibrations of this kind exist for other potentials, such as the anisotropy energy, or for skyrmions of higher degree.

\paragraph{Outline.}
In Section~\ref{sec:profile} we introduce Bloch fields and the profile equation \eqref{profileode} and fix a solution $f$; 
its existence is proved in Appendix~\ref{ex-existence} by a continuation argument in $\mu$, 
starting from a result of Ibrahim and Shimizu~\cite{IS}. 
In Sections~\ref{sec:barriers} and~\ref{sec:target}, we prove $a_0>1$ and the triangle inequalities, 
and construct $H$ and $F$. 
Section~\ref{sec:calibration} contains the pointwise identity; Section~\ref{sec:identity} its integration over $\R^2$; 
and Section~\ref{sec:equality} the equality case and the exact range of the method (Theorem~\ref{conditional}, Proposition~\ref{sharpness}).
Section~\ref{sec:branch} determines the calibration endpoint $\mu_{\rm cal}$ of Theorem~\ref{rangetheorem}.
In Section~\ref{sec:original} we return to the original parameters and complete the proofs of Theorems~\ref{originaltheorem} and~\ref{rangetheorem}.

\section{Bloch fields and the profile equation}
\label{sec:profile}

For $\mu>0$, write $E_\mu=E_{\mu,1}$. The Bloch convention is
\[
 B_f(r,\theta)=\sin f(r)e_\theta+\cos f(r)e_3,
 \qquad e_\theta=(-\sin\theta,\cos\theta,0).
\]
Using \eqref{dmiext} gives
\[
 E_\mu(B_f)=\pi\int_0^\infty\left[
  rf'^2+\frac{\sin^2f}{r}+4\mu r\sin^2f f'
                    +2r(1-\cos f)\right]\dd r.
\]
Its Euler equation and desired profile class are
\begin{equation}
 \begin{gathered}
 -f''-f'/r+\sin f\cos f/r^2-2\mu\sin^2f/r+\sin f=0,\\
 f(0)=\pi,\quad f(\infty)=0,\quad0<f<\pi,\quad f'<0,\\
 \int_0^\infty(f'^2+\sin^2f/r^2)r\dd r<\infty,
 \qquad\int_0^\infty(1-\cos f)r\dd r<\infty.
 \end{gathered}                                      \label{profileode}
\end{equation}
This derivation fixes the DMI sign. Its charge is
$\frac12\int_0^\infty f'\sin f\dd r=-1$.

\begin{proposition}[Existence]
\label{profileexists}
For each $0<\mu\le\sqrt2$ there exists at least one profile in
\eqref{profileode}, regular at the origin.
\end{proposition}

The proof is given in Appendix~\ref{ex-existence}. 
It is a continuation argument in $\mu$. 
A priori bounds and a compactness argument show that the set of $\mu\in(0,\sqrt2]$ for which such a profile exists is closed in $(0,\sqrt{2}]$. 
A positive supersolution of the linearized equation shows that every profile is nondegenerate, so the implicit function theorem makes this set open. 
Finally, Theorem~A of Ibrahim--Shimizu~\cite{IS}, with their $\beta=1$, coupling $\mathfrak r=\sqrt2\mu$, and radius $\rho=r/\sqrt2$, shows that the set contains $\mu=1/2$. 
Since $(0,\sqrt2]$ is connected, the set is all of $(0,\sqrt2]$.
The appendix does not use the calibration, minimality, or symmetry results proved later. 
From now on we fix one such profile.
The calibration of Sections~\ref{sec:barriers}--\ref{sec:equality} is built from a given profile. 
It requires $a_0\ge1$ (and $\mu>3/4$ if $a_0=1$); this is proved for $\mu\le1$ in Lemma~\ref{originslope} and imposed as a hypothesis otherwise.

\begin{lemma}[Asymptotics at the origin and at infinity]\label{profileendpoints}
Every profile in \eqref{profileode} has a finite slope $a_0=-f'(0)>0$,
is the unique regular-origin solution for that datum, and satisfies
\begin{equation}
 f(r)=\pi-a_0r+c_3r^3+O(r^5),\qquad
 c_3=\frac{a_0+(2/3)a_0^3-2\mu a_0^2}{8}.       \label{originseries}
\end{equation}
The remainder can be differentiated once, 
and the remainder and its derivative are bounded uniformly for $(\mu,a_0)$ in compact subsets of $(0,\infty)\times(0,\infty)$.
At infinity, there is $c_f>0$ such that
\begin{equation}
 f=c_fK_1+o(K_1/r^2),\qquad f'=c_fK_1'+o(K_1/r^2),
                                                        \label{masstail}
\end{equation}
where $K_1$ is the modified Bessel function of the second kind of order one.
\end{lemma}
\begin{proof}
The class \eqref{profileode} is the class $\mathcal P_\mu$ of \eqref{ex-class}. 
The assertions \eqref{originseries} at the origin are proved in Lemma~\ref{ex-origin-lemma}, and \eqref{masstail} in Lemma~\ref{ex-tail-lemma}; 
both proofs use only the differential equation.
\end{proof}

Up to a constant factor, $K_1$ is the only solution of
\[
-y''-\frac{y'}{r}+\frac{y}{r^2}+y=0,
\]
the linearization of the equation in \eqref{profileode} at $f=0$, that decays at infinity. 
We will use the expansions
\begin{equation}\label{besselprecision}
K_1(r)=\sqrt{\frac{\pi}{2r}}\,e^{-r}\Big(1+\frac{3}{8r}+O(r^{-2})\Big),\qquad
-\frac{K_1'}{K_1}=1+\frac{1}{2r}+O(r^{-2}),
\end{equation}
which follow, for example, from the integral representation $K_1(r)=\int_0^\infty e^{-r\cosh t}\cosh t\,dt$.

\section{The slope at the origin and the triangle inequalities}
\label{sec:barriers}

The next lemma bounds the slope of a profile at the origin from below. 
It shows that the hypothesis $a_0\ge1$ of Lemma~\ref{triangles} and of Sections~\ref{sec:target}--\ref{sec:equality} is satisfied for $\mu\le1$; 
it is also used in Section~\ref{sec:branch}.

\begin{lemma}\label{originslope}
Let $0<\mu \leq 1$ and let $f$ be a profile in \eqref{profileode}. Then $a_0=-f'(0)>1$.
\end{lemma}
\begin{proof}
For $\psi\in C^2(0,\infty)$ let
\[
\mathcal L_\mu\psi=-\psi''-\frac{\psi'}{r}+\frac{\sin\psi\cos\psi}{r^2}
-\frac{2\mu\sin^2\psi}{r}+\sin\psi
\]
be the left side of the equation in \eqref{profileode}, 
so that $\mathcal L_\mu f=0$. 
We construct strict supersolutions $\psi_\lambda$ with slope $-\lambda$ at the origin. 
If $a_0\le1$, they can be slid down to $\lambda=a_0$ while staying above $f$, 
and this contradicts the expansion of $f$ at the origin.

\emph{Supersolutions.} Let $\chi=(1+\sqrt3)/16$ and
\[
\psi_*(r)=2\arctan\frac{2}{r(1+\chi r^2)},
\]
a modification of the Belavin--Polyakov profile $2\arctan(2/r)$. 
It decreases from $\pi$ to $0$, satisfies $\psi_*(r)=\pi-r+(1/12-\chi)r^3+O(r^5)$ as $r\to0$, 
and $\psi_*(r)\sim4/(\chi r^3)$ as $r\to\infty$. 
A direct computation using $\tan(\psi_*/2)=2/(r(1+\chi r^2))$ gives
\begin{equation}\label{supfactor}
\frac{\mathcal L_1\psi_*}{\sin\psi_*}
=\frac{(\chi r^2-8\chi+1)(\chi r^3+r-2)(\chi r^3+r+2)}
{(1+\chi r^2)(\chi^2r^6+2\chi r^4+r^2+4)}\ge0 .
\end{equation}
Indeed, the first two factors in the numerator are increasing and both vanish at $r=2(\sqrt3-1)$ (this is how $\chi$ was chosen), 
so their product is nonnegative; the other factors are positive. 
Moreover,
\[
\frac{\sin\psi_*}{r}=\frac{4(1+\chi r^2)}{r^2(1+\chi r^2)^2+4}\le1,
\]
since the denominator minus the numerator is $r^2\big((1+\chi r^2)^2-4\chi\big)>0$.

For $0<\lambda<1$, put $\psi_\lambda(r)=\psi_*(\lambda r)$. 
With $x=\lambda r$, a direct computation gives
\[
\mathcal L_\mu\psi_\lambda(r)=\lambda^2\mathcal L_1\psi_*(x)
+\frac{2\lambda(\lambda-\mu)\sin^2\psi_*(x)}{x}+(1-\lambda^2)\sin\psi_*(x).
\]
The first term is nonnegative by \eqref{supfactor}. 
If $\lambda\ge\mu$, so is the second, and the third is positive. 
If $\lambda<\mu$, the bound $\sin\psi_*(x)/x\le1$ shows that the last two terms are at least
$(1+\lambda^2-2\lambda\mu)\sin\psi_*(x)\ge(1-\lambda)^2\sin\psi_*(x)>0$. 
Hence $\mathcal L_\mu\psi_\lambda>0$ on $(0,\infty)$ for every $0<\lambda<1$.

\emph{Sliding.} Let $\lambda_0=\min\{a_0,1\}$. 
We show that $\psi_\lambda>f$ on $(0,\infty)$ for every $\lambda\in(0,\lambda_0)$. 
By \eqref{originseries} and the expansion of $\psi_*$,
\[
\psi_\lambda(r)-f(r)=(a_0-\lambda)r+O(r^3)\qquad(r\to0),
\]
uniformly in $\lambda\in(0,1]$. 
Since $\psi_*$ is decreasing, $\psi_\lambda\ge\psi_1$ for $\lambda\le1$; 
moreover, $\psi_1$ decays like $r^{-3}$, whereas $f$ decays exponentially by \eqref{masstail}. 
Fix $\lambda_2\in(0,\lambda_0)$. 
Since $a_0-\lambda\ge a_0-\lambda_2>0$ for $\lambda\le\lambda_2$, 
there are $0<r_1<R_1$ such that $\psi_\lambda>f$ on $(0,r_1]\cup[R_1,\infty)$ for all $\lambda\in(0,\lambda_2]$.

Let $\Lambda$ be the set of $\lambda\in(0,\lambda_2]$ such that $\psi_\lambda>f$ on $[r_1,R_1]$. 
It contains all small $\lambda$, 
because $\psi_\lambda\to\pi$ uniformly on $[r_1,R_1]$ as $\lambda\to0$ while $f<\pi$. 
It is open in $(0,\lambda_2]$ because $\psi_\lambda-f$ is continuous in $(\lambda,r)$ and $[r_1,R_1]$ is compact. 
It is also closed in $(0,\lambda_2]$. 
Indeed, let $\lambda_n\in\Lambda$ converge to
$\bar\lambda\in(0,\lambda_2]$. 
Then $w=\psi_{\bar\lambda}-f\ge0$ on $(0,\infty)$. If $w(r_*)=0$ for some $r_*>0$, 
then $w$ has a minimum at $r_*$, so $w'(r_*)=0$ and $w''(r_*)\ge0$. 
Since $\psi_{\bar\lambda}$ and $f$ have the same value and the same derivative at $r_*$, 
and only the term $-\psi''$ of $\mathcal L_\mu\psi$ involves $\psi''$,
\[
\mathcal L_\mu\psi_{\bar\lambda}(r_*)
=\mathcal L_\mu\psi_{\bar\lambda}(r_*)-\mathcal L_\mu f(r_*)=-w''(r_*)\le0,
\]
which contradicts $\mathcal L_\mu\psi_{\bar\lambda}>0$. 
Hence $w>0$ and $\bar\lambda\in\Lambda$. 
Since $(0,\lambda_2]$ is connected, $\Lambda=(0,\lambda_2]$; 
as $\lambda_2<\lambda_0$ was arbitrary, the claim follows.

\emph{Conclusion.} Suppose $a_0\le1$. 
Then $\lambda_0=a_0$, and letting $\lambda\uparrow a_0$ gives $\psi_{a_0}\ge f$. 
But $\psi_{a_0}$ and $f$ have the same slope at the origin, 
and by \eqref{originseries} and the expansion of $\psi_*$,
\[
\psi_{a_0}(r)-f(r)=-a_0\Big(\chi a_0^2-\frac{\mu a_0}{4}+\frac18\Big)r^3+O(r^5).
\]
Since $\mu \leq 1$ and $\chi>1/8$, the bracket is at least $\chi a_0^2-a_0/4+1/8=(a_0-1)^2/8+(\chi-1/8)a_0^2>0$. 
Hence $\psi_{a_0}<f$ for small $r>0$, a contradiction. 
Therefore $a_0>1$.
\end{proof}

Let $f$ be a profile in \eqref{profileode}. 
As functions of $r>0$, put
\[
A=-f',\qquad b=\frac{\sin f}{r},\qquad M=2\sin\frac f2,
\]
and write $u=\sin(f/2)$, $v=\cos(f/2)$, $c=\cos f$. 
Since $0<f<\pi$ and $f'<0$, the functions $A$, $b$, $M$, $u$, and $v$ are positive. 
The next lemma states that $A$, $b$, and $M$ are the side lengths of a nondegenerate triangle,
that is, $|M^2-A^2-b^2|<2Ab$.
In Section~\ref{sec:target} this becomes the bound $|H|<1$. Set
\[
W_0=A+b-M,\qquad W_1=M-A+b,\qquad W_2=M+A-b .
\]
Differentiating and using the equation in \eqref{profileode} to eliminate $f''$, we obtain
\begin{align}
W_0'&=\Big(v-\frac{1+c}{r}\Big)W_0+\frac{2uv^2}{r}\,(4\mu u-3),\label{Wzero}\\
W_1'&=\Big(v-\frac{1-c}{r}\Big)W_1+\frac{2u}{r}\big(3u^2-1-4\mu u(1-u^2)\big),\label{Wone}\\
W_2'&=\Big(-v-\frac{1-c}{r}\Big)W_2+\frac{2u}{r}\big(3u^2-1+4\mu u(1-u^2)\big).\label{Wtwo}
\end{align}

\begin{lemma}[Strict triangle inequalities]\label{triangles}
Let $\mu>0$ and let $f$ be a profile in \eqref{profileode} with $a_0=-f'(0)\ge1$; if $a_0=1$, assume in addition that $\mu>3/4$.
Then $|A-b|<M<A+b$ for all $r>0$; 
equivalently, $W_0$, $W_1$ and $W_2$ are positive.
\end{lemma}
\begin{proof}
Write \eqref{Wzero}--\eqref{Wtwo} as $W_i'=d_iW_i+g_i$ ($i=0,1,2$), 
and put the integrating factor $P_i=\exp(-\int_1^r d_i)>0$ so that $(P_iW_i)'=P_ig_i$. 
We use the following observation, in which we drop the index $i$.

\smallskip\noindent
$(\ast)$ \emph{Assume that $PW\to\ell_0\ge0$ as $r\to0$ and $PW\to0$ as $r\to\infty$, 
that $Pg$ is integrable on $(0,\infty)$, 
and that for some $r_1\ge0$ we have $Pg>0$ on $(0,r_1)$ and $Pg<0$ on $(r_1,\infty)$. 
Then $W>0$ on $(0,\infty)$.}

\smallskip\noindent
Indeed, $PW(r)=\ell_0+\int_0^rPg>0$ for $0<r\le r_1$, and $PW(r)=-\int_r^\infty Pg>0$ for $r>r_1$. 
In words: $PW$ starts at $\ell_0\ge0$, increases while the source is positive, 
and then decreases to $0$ at infinity, 
so it stays positive in between.

\emph{Signs of the sources.} Since $f$ decreases from $\pi$ to $0$, $u$ decreases strictly from $1$ to $0$. 
The source in \eqref{Wzero} has the sign of $4\mu u-3$: it is negative everywhere if $\mu\le3/4$, and if $\mu>3/4$, 
it is positive for $u>3/(4\mu)$ and negative for $u<3/(4\mu)$. 
The sources in \eqref{Wone} and \eqref{Wtwo} have the signs of $p_-(u)$ and $p_+(u)$, 
respectively, where $p_\mp(u)=3u^2-1\mp4\mu u(1-u^2)$; note that $p_\mp(0)=-1$, $p_\mp(1)=2$. 
Since $p_-''(u)=6+24\mu u>0$, the polynomial $p_-$ is convex and has exactly one zero in $(0,1)$. 
Since $p_+'(u)=6u+4\mu-12\mu u^2$ is positive at $u=0$ and has exactly one positive zero, 
$p_+$ increases and then possibly decreases on $[0,1]$; as $p_+(1)>0$, it also has exactly one zero in $(0,1)$. 
Thus, as $r$ increases, each source changes sign at most once, from positive to negative.

\emph{Behavior at the origin.} By Lemma~\ref{profileendpoints}, $A,b\to a_0$, $M\to2$, $v\to0$, and $c\to-1$ as $r\to0$. 
Hence $W_0\to2a_0-2\ge0$, since $a_0\ge1$, and $W_1,W_2\to2$.
The coefficient $v-(1+c)/r$ in \eqref{Wzero} is bounded near $0$, 
while the coefficients in \eqref{Wone} and \eqref{Wtwo} are $-2/r+O(r)$. 
So the integrating factors satisfy $P_0\to p_0>0$ and $P_1,P_2\sim p_ir^2$ with $p_i>0$. 
Consequently $P_0W_0\to p_0(2a_0-2)\geq 0$, $P_1W_1,P_2W_2\to0$, and all $P_ig_i$ are $O(r)$.

\emph{Behavior at infinity.} By Lemma~\ref{profileendpoints}, $f$ and $f'$ are $O(K_1)$. 
Hence $u$, $W_0$, $W_1$ and $W_2$ are $O(K_1)=O(r^{-1/2}e^{-r})$, while $1-v$ and $1-c$ are $O(K_1^2)$. 
The coefficients in \eqref{Wzero}--\eqref{Wtwo} are therefore $1-2/r$, $1$ and $-1$ up to integrable errors, 
so $P_0=O(r^2e^{-r})$, $P_1=O(e^{-r})$ and $P_2=O(e^r)$. 
It follows that every $P_iW_i$ tends to $0$; for $i=2$, the growth of $P_2$ is compensated since $P_2W_2=O(r^{-1/2})$. 
Moreover, $g_i=O(K_1/r)$, so $P_0g_0$ and $P_1g_1$ decay exponentially and $P_2g_2=O(r^{-3/2})$. 
In particular, each $P_ig_i$ is integrable on $(0,\infty)$.

Now $(\ast)$ applies to $W_1$ and $W_2$ with $\ell_0=0$, 
and to $W_0$ with $\ell_0=p_0(2a_0-2)\ge0$: if $a_0>1$, then $\ell_0>0$, and $r_1=0$ is allowed; 
if $a_0=1$, then $\ell_0=0$, and $r_1>0$ because $\mu>3/4$. 
Hence $W_0,W_1,W_2>0$.
\end{proof}

\section{The functions \texorpdfstring{$H$}{H} and \texorpdfstring{$F$}{F}}
\label{sec:target}

In this section, we use the profile $f$ to define two functions $H$ and $F$ of the polar angle $t\in[0,\pi]$ on the target sphere, 
where $m_3=\cos t$. 
In Section~\ref{sec:calibration}, 
they are the coefficients of a pointwise identity for the energy density of an arbitrary map.
In this section and the next three, we assume that $a_0\ge1$, and that $\mu>3/4$ if $a_0=1$; 
by Lemma~\ref{originslope}, this holds for every profile when $\mu\le1$.

Since $f$ decreases strictly from $\pi$ to $0$, it has an inverse $R(t)=f^{-1}(t)$, $0<t<\pi$; 
the south pole $t=\pi$ corresponds to $r=0$ and the north pole $t=0$ to $r=\infty$. 
Set
\[
s=\sin t,\quad c=\cos t,\quad V=1-c,\quad A(t)=-f'(R(t)),\quad
\beta(t)=\frac{s}{R(t)},\quad b_0(t)=\frac1{R(t)}.
\]
Thus $A$ and $\beta$ are the lengths of $\partial_rB_f$ and $r^{-1}\partial_\theta B_f$ at the radius where $B_f$ has polar angle $t$. Define
\begin{equation}\label{targetHF}
H=\frac{2V-A^2-\beta^2}{2A\beta},\qquad
F=-s(\beta+HA)=\frac{A^2-\beta^2-2V}{2b_0}.
\end{equation}
These formulas are chosen so that the identity of Section~\ref{sec:calibration} is an equality for the Bloch field $B_f$. 
That identity also requires the following compatibility condition. 
Differentiating $f(R(t))=t$ gives $R_t=-1/A$, so $(b_0)_t=b_0^2/A$; 
moreover, $AA_t=f''(R)$, and the equation in \eqref{profileode} at $r=R(t)$ gives
\[
AA_t=Ab_0+scb_0^2-2\mu s^2b_0+s .
\]
Differentiating the last expression in \eqref{targetHF} and using these two relations, we obtain
\begin{equation}\label{compatibility}
F_t=A+H\beta-2\mu s^2 .
\end{equation}
Only the profile $f$ enters here; 
in Section~\ref{sec:calibration} the functions $H$ and $F$ are evaluated along an arbitrary map, 
which need not satisfy any equation.

\emph{The interior.} For $0<t<\pi$, Lemma~\ref{triangles} at $r=R(t)$ gives $|A-\beta|<2\sin(t/2)<A+\beta$. 
Squaring and using $2V=4\sin^2(t/2)$, we get $|2V-A^2-\beta^2|<2A\beta$, that is, $-1<H<1$.

\emph{The south pole.} As $t\to\pi$, $R\to0$. 
By \eqref{originseries}, $A$ and $\beta$ are $a_0+O(R^2)$, and $R=(\pi-t)/a_0+O((\pi-t)^3)$. 
Since $2V\to4$, it follows that $H\to2/a_0^2-1$ and $F=-2(\pi-t)/a_0+O((\pi-t)^3)$.

\emph{The north pole.} As $t\to0$, $R\to\infty$, and $t=f(R)=O(e^{-R})$. 
By \eqref{masstail} and \eqref{besselprecision},
\[
\frac At=1+\frac1{2R}+O(R^{-2}),\qquad \frac\beta t=\frac1R\big(1+O(t^2)\big),\qquad
\frac{2V}{t^2}=1+O(t^2).
\]
The leading terms of $2V$ and $A^2$ cancel, so $2V-A^2-\beta^2=-t^2/R+O(t^2/R^2)+O(t^4)$, 
while $2A\beta=(2t^2/R)(1+O(R^{-1}))$. 
Since $Rt^2\to0$, the quotient $H$ tends to $-1/2$; this is why the second term in \eqref{besselprecision} is needed. 
In the same way, $F=R(A^2-\beta^2-2V)/2=\frac12t^2(1+O(R^{-1}))$. In summary,
\begin{equation}\label{polelimits}
H(0)=-\frac12,\qquad H(\pi)=\frac2{a_0^2}-1,\qquad
F=\frac12t^2\big(1+O(R^{-1})\big)\quad(t\to0).
\end{equation}
 
Since $a_0\ge1$, we have $H(\pi)=2/a_0^2-1\in(-1,1]$, with $H(\pi)=1$ exactly when $a_0=1$. 
Together with the bound in the interior, this gives
\begin{equation}\label{calgap}
H\in C([0,\pi]),\qquad |H|<1\ \text{on }(0,\pi),\qquad
\delta_\mu=1-\max_{[0,\pi]}|H|\ \begin{cases}>0&\text{if }a_0>1,\\ =0&\text{if }a_0=1.\end{cases}
\end{equation}
Finally, $F$ vanishes at both poles; $F_t$ and $F/s$ are bounded on $(0,\pi)$ and are $O(s)$ as $t\to0$; 
and $A$ and $\beta$ are bounded and are $O(s)$ as $t\to0$. 
These bounds follow from \eqref{compatibility} and the asymptotics above, and they are used in
Section~\ref{sec:identity}. 
All constants depend on the fixed profile $f$; they are not uniform as $\mu\to0$, and $\delta_\mu\to0$ as $a_0\to1$.

\section{A pointwise identity for the energy density}
\label{sec:calibration}

Let $m$ be a smooth map. Where $m\ne\pm e_3$, write locally
\[
m=(\sin t\,n,\cos t),\qquad n=(\cos\phi,\sin\phi),
\]
and set
\[
u=\nabla t,\qquad v=\sin t\,\nabla\phi,\qquad X=u-A(t)\,Jn,\qquad Y=v-\beta(t)\,n,
\]
with $A$ and $\beta$ as in Section~\ref{sec:target}. 
The vectors $\tau_1=(\cos t\,n,-\sin t)$ and $\tau_2=(Jn,0)$ form an orthonormal basis of the tangent plane of $\mathbb S^2$ at $m$, 
with $\tau_1\times\tau_2=m$ and $\partial_im=u_i\tau_1+v_i\tau_2$. Hence
\[
|\nabla m|^2=|u|^2+|v|^2,\qquad m\cdot(\partial_1m\times\partial_2m)=\det(u,v).
\]
The Bloch field $B_f$ has $u=AJn$ and $v=\beta n$; that is, $X=Y=0$. 
Thus $X$ and $Y$ measure the deviation of $\nabla m$ from the gradient of $B_f$ at the same polar angle.

With the DMI term written as in \eqref{dmiext}, the energy density of $m$ is
\[
e_{\rm cal}=\frac12\big(|u|^2+|v|^2\big)-2\mu s^2\,Jn\cdot u+V,
\]
with $s=\sin t$ and $V=1-\cos t$ as in Section~\ref{sec:target}; its integral is $E_\mu(m)$. By the chain rule,
\begin{equation}\label{fluxlocal}
\operatorname{div}(F\,Jn)=F_t\,Jn\cdot u-\frac Fs\,n\cdot v .
\end{equation}

\begin{lemma}[Calibration identity]\label{lem:calibration}
At every point where $m\ne\pm e_3$,
\begin{equation}\label{CAL}
e_{\rm cal}=H\det(u,v)+\operatorname{div}(F\,Jn)+\mathcal R_\mu(X,Y),\qquad
\mathcal R_\mu(X,Y)=\frac12\big(|X|^2+|Y|^2\big)-H\det(X,Y),
\end{equation}
where $H$ and $F$ are evaluated at $t$. Moreover,
\begin{equation}\label{remainderpositive}
\mathcal R_\mu(X,Y)\ge\frac{1-|H(t)|}{2}\big(|X|^2+|Y|^2\big)\ge\frac{\delta_\mu}{2}\big(|X|^2+|Y|^2\big).
\end{equation}
\end{lemma}

\begin{proof}
Fix a point, and hence $t$ and $n$. By \eqref{fluxlocal},
\[
q(u,v)=e_{\rm cal}-H\det(u,v)-\operatorname{div}(F\,Jn)
\]
is a quadratic polynomial in $(u,v)\in\mathbb R^2\times\mathbb R^2$ with a quadratic part $\frac12(|u|^2+|v|^2)-H\det(u,v)$. 
A quadratic polynomial whose value and gradient vanish at a point equals its quadratic part evaluated at the displacement from that point. 
The displacement from the Bloch point $(AJn,\beta n)$ is $(X,Y)$, so it suffices to show that $q$ and its gradient vanish there.

Since $\nabla_u\det(u,v)=-Jv$, $\nabla_v\det(u,v)=Ju$, and $J^2=-1$, at the Bloch point
\[
\nabla_uq=\big(A+H\beta-2\mu s^2-F_t\big)Jn,\qquad
\nabla_vq=\Big(\beta+HA+\frac Fs\Big)n .
\]
These vanish by \eqref{compatibility} and by the definition of $F$ in \eqref{targetHF}. 
Using $\det(Jn,n)=-1$ and these two relations, the value of $q$ at the Bloch point is
\[
\frac12(A^2+\beta^2)+V+HA\beta-(A+H\beta)A-(\beta+HA)\beta
=V-\frac12(A^2+\beta^2)-HA\beta,
\]
which vanishes by the definition of $H$ in \eqref{targetHF}. This proves
\eqref{CAL}.

Finally, $|\det(X,Y)|\le|X||Y|\le\frac12(|X|^2+|Y|^2)$ gives the first inequality in \eqref{remainderpositive}, 
and $|H|\le1-\delta_\mu$ by \eqref{calgap} gives the second.
\end{proof}

The identity \eqref{CAL} is algebraic: it holds for every smooth map wherever $m\ne\pm e_3$ and uses no property of $m$. 
For the exchange energy alone, the Bogomol'nyi argument~\cite{BP75,Bog76} gives an identity of the same type, 
$\frac12(|u|^2+|v|^2)=-\det(u,v)+\frac12|u-Jv|^2$, with $H\equiv-1$ and a remainder that is only positive semidefinite. 
Here $H$ depends on $t$ and $|H(t)|<1$ for $0<t<\pi$, so the remainder is positive definite at every point where $m\ne\pm e_3$. 
For $B_f$, the topological density is $\det(AJn,\beta n)=-A\beta<0$, in agreement with its degree $-1$. 

\section{Integration of the identity}
\label{sec:identity}

In this section we integrate \eqref{CAL} over $\R^2$ for every $m\in\Aminus$. 
Since such maps need not be smooth, we first define $t$, $u$ and $v$ for them. 
Lemma~\ref{weighteddegree} computes the integral of the first term in \eqref{CAL}, 
and Lemma~\ref{globalidentity} shows that the divergence term integrates to zero.

Let $m\colon\R^2\to\mathbb S^2$ with $m-e_3\in H^1$, and write $m=(w,z)$ with $w\in\R^2$. 
The polar angle $t=\arccos z$ is the geodesic distance from $m$ to $e_3$ on $\mathbb S^2$. 
Since geodesic distance on $\mathbb S^2$ is at most $\pi/2$ times Euclidean distance in $\R^3$, 
$t$ is a Lipschitz function of $m\in\mathbb S^2$ with constant $\pi/2$, and $t\le(\pi/2)|m-e_3|$. 
Extending it to a Lipschitz function on $\R^3$ and using the chain rule for Sobolev maps, we obtain $t\in H^1$; set $u=\nabla t$.

On $\Omega=\{|z|<1\}=\{m\ne\pm e_3\}$ put $n=w/|w|$ and $v_i=(Jn)\cdot\partial_iw$. 
Where a local azimuth $\phi$ with $n=(\cos\phi,\sin\phi)$ exists, 
we have $\partial_iw=\cos t\,\partial_it\,n+\sin t\,\partial_i\phi\,Jn$, so $v=\sin t\,\nabla\phi$ as in Section~\ref{sec:calibration}; 
the definition of $v$ does not require a global azimuth. 
Since $|m|=1$, the derivatives $\partial_im$ are tangent to $\mathbb S^2$ almost everywhere, 
and on $\Omega$ they satisfy $\partial_im=u_i\tau_1+v_i\tau_2$ with $\tau_1$, $\tau_2$ as in Section~\ref{sec:calibration}. 
Hence, almost everywhere on $\Omega$,
\[
|\nabla m|^2=|u|^2+|v|^2,\qquad j_m:=m\cdot(\partial_1m\times\partial_2m)=\det(u,v).
\]
On the pole sets $\{m=\pm e_3\}$, all weak derivatives of $m$ vanish almost everywhere 
because the gradient of a Sobolev function vanishes almost everywhere on each of its level sets. 
There $j_m=0$, and the energy density reduces to the Zeeman density $\frac12|m-e_3|^2$, 
which is $0$ on $\{m=e_3\}$ and $2$ on $\{m=-e_3\}$. 
This produces the term $2|\{m=-e_3\}|$ in \eqref{globalCAL} below.

For smooth $m$, the next lemma is the change of variables formula: 
$j_m\,dx$ is the pullback of the area form $\sin t\,dt\wedge d\phi$ of $\mathbb S^2$, 
so $\int G(t)\,j_m$ is $Q(m)$ times $\int_{\mathbb S^2}G=2\pi\int_0^\pi G(t)\sin t\dd t$. 
For $G\equiv1$, it is the definition of $Q(m)$.

\begin{lemma}[Weighted degree]\label{weighteddegree}
Let $m\colon\R^2\to\mathbb S^2$ with $m-e_3\in H^1$, and let $G$ be continuous on $[0,\pi]$. Then
\[
\int G(t)\,j_m=2\pi\,Q(m)\int_0^\pi G(t)\sin t\dd t .
\]
\end{lemma}
\begin{proof}
Let $\bar G=\frac12\int_0^\pi G(t)\sin t\dd t$ be the mean of $G$ over $\mathbb S^2$. 
Since $\int j_m=4\pi Q(m)$, the claim holds for the constant $\bar G$, so it suffices to show that $\int(G(t)-\bar G)\,j_m=0$. Put
\[
K(t)=\int_0^t\big(G(q)-\bar G\big)\sin q\dd q .
\]
Then $K(0)=K(\pi)=0$ and $|K(t)|\le2\|G\|_\infty\sin^2t$. 
Away from the poles, $d(K(t)\,d\phi)=(G(t)-\bar G)\sin t\,dt\wedge d\phi$, so $K\,d\phi$ is a primitive of $G-\bar G$ times the area form. 
Since $d\phi$ pulls back to $(w_1\nabla w_2-w_2\nabla w_1)/|w|^2$ and $|w|^2=1-z^2$, the pullback of $K\,d\phi$ by $m$ is
\[
\alpha=g(z)\,(w_1\nabla w_2-w_2\nabla w_1),\qquad g(z)=\frac{K(\arccos z)}{1-z^2}.
\]
The bound on $K$ gives $|g|\le2\|G\|_\infty$, so $\alpha$ is defined on all of $\R^2$ (it vanishes on the pole sets), 
and $\alpha\in L^1$ because $w,\nabla w\in L^2$. We claim that
\begin{equation}\label{primitivecurl}
\partial_1\alpha_2-\partial_2\alpha_1=\big(G(t)-\bar G\big)\,j_m\qquad\text{in }\mathcal
D'(\R^2).
\end{equation}
For smooth $m$, this says that pullback commutes with $d$: locally $\alpha=K(t)\nabla\phi$, 
so $\partial_1\alpha_2-\partial_2\alpha_1=K'(t)(\partial_1t\,\partial_2\phi-\partial_2t\, \partial_1\phi)=(G(t)-\bar G)\,j_m$.

\emph{Polynomial $G$.} Suppose first that $G(t)=p(\cos t)$ for a polynomial $p$. 
Then $K(\arccos z)=\int_z^1(p(q)-\bar G)\dd q$ is a polynomial in $z$ that vanishes at $z=\pm1$, so $g$ is a polynomial, 
and $\alpha$ is a polynomial expression in $m$ and $\nabla m$ that makes sense for every $\R^3$-valued map. 
For a smooth $\R^3$-valued map $m=(w_1,w_2,z)$, put $\omega_i=w_1\partial_iw_2-w_2\partial_iw_1$ 
and $D=\partial_1w_1\partial_2w_2-\partial_1w_2\partial_2w_1$. 
Since $\partial_1\partial_2=\partial_2\partial_1$, the second derivatives cancel and
\[
\partial_1\alpha_2-\partial_2\alpha_1=g'(z)\big(\partial_1z\,\omega_2-\partial_2z\,\omega_1\big)+2g(z)\,D,
\]
a polynomial in $m$ times products of two first derivatives. 
Expanding $j_m=m\cdot(\partial_1m\times\partial_2m)$ gives $\partial_1z\,\omega_2-\partial_2z\,\omega_1=zD-j_m$ for every map. 
At points where $|m|=1$ and $m\cdot\partial_im=0$, 
the vector $\partial_1m\times\partial_2m$ is orthogonal to the tangent plane, hence equals $j_m\,m$; 
its third component gives $D=z\,j_m$.
Therefore, at such points,
\[
\partial_1\alpha_2-\partial_2\alpha_1=\big(2zg(z)-(1-z^2)g'(z)\big)\,j_m=\big(p(z)-\bar G\big)\,j_m,
\]
where the last equality follows from differentiating $(1-z^2)g(z)=\int_z^1(p-\bar G)\dd q$. 
Now let $m_\varepsilon=m*\rho_\varepsilon$ with a nonnegative mollifier $\rho_\varepsilon$. 
Then $|m_\varepsilon|\le1$, $m_\varepsilon-e_3\to m-e_3$ in $H^1$, 
and $m_\varepsilon\to m$ almost everywhere along a subsequence. 
Hence $\alpha(m_\varepsilon)\to\alpha(m)$ in $L^1$, and the curl of $\alpha(m_\varepsilon)$, 
computed for $\R^3$-valued maps, 
converges in $L^1$ to the same expression evaluated at $m$: its coefficients are bounded and converge almost everywhere, 
and the products of first derivatives converge in $L^1$. 
Since $|m|=1$ and $m\cdot\partial_im=0$ almost everywhere, the limit is $(p(z)-\bar G)\,j_m$. 
This proves \eqref{primitivecurl} for polynomial $G$.

\emph{Continuous $G$.} Choose polynomials $p_k$ with $p_k(z)\to G(\arccos z)$ uniformly on $[-1,1]$, 
and put $G_k(t)=p_k(\cos t)$. 
Then $\bar G_k\to\bar G$, and the bound on $K$, applied to $G_k-G$, shows that $g_k\to g$ uniformly. 
Hence, the corresponding forms converge in $L^1$, and $(G_k(t)-\bar G_k)\,j_m\to(G(t)-\bar G)\,j_m$ in $L^1$ because $|j_m|\le\frac12|\nabla m|^2$. 
This proves \eqref{primitivecurl} for continuous $G$.

\emph{Conclusion.} Let $\chi_R$ be a cutoff function with $\chi_R=1$ on $B_R$, $\operatorname{supp}\chi_R\subset B_{2R}$, 
and $|\nabla\chi_R|\le C/R$. By \eqref{primitivecurl},
\[
\Big|\int\chi_R\big(G(t)-\bar G\big)\,j_m\Big|
=\Big|\int\big(\alpha_1\partial_2\chi_R-\alpha_2\partial_1\chi_R\big)\Big|
\le\frac CR\,\|\alpha\|_{L^1}\longrightarrow0 .
\]
By the dominated convergence theorem, the left side tends to $\big|\int(G(t)-\bar G)\,j_m\big|$, which is therefore zero.
\end{proof}

\begin{lemma}[Energy identity]\label{globalidentity}
Let $m\colon\R^2\to\mathbb S^2$ with $m-e_3\in H^1$, and let
\[
T=F(t)\,Jn=\frac{F(t)}{\sin t}\,Jw\quad\text{on }\Omega,\qquad T=0\quad\text{on the pole sets}.
\]
Then $T\in L^2$, $\diver T\in L^1$, and $\int\diver T=0$. If moreover $m\in\Aminus$, then
\begin{equation}
E_\mu(m)=C_\mu+\int_\Omega\mathcal R_\mu(X,Y)\dd x+2|\{m=-e_3\}|,
\qquad C_\mu=-2\pi\int_0^\pi H(t)\sin t\dd t.    \label{globalCAL}
\end{equation}
\end{lemma}
\begin{proof}
Since $|m-e_3|^2=2(1-\cos t)\ge\sin^2t$, we have $\sin t\le|m-e_3|$; moreover, $|m-e_3|$ is bounded below on any set where $t$ is bounded away from $0$. 
Hence, every function $g$ that is bounded on $(0,\pi)$ and is $O(\sin t)$ as $t\to0$ satisfies $|g(t)|\le C|m-e_3|$ on $\mathbb S^2\setminus\{\pm e_3\}$. 
By Section~\ref{sec:target}, this applies to $F_t$, $F/s$, $A$, and $\beta$.

\emph{$T$ is a Lipschitz function of $m$.} Let $\Phi(m)=F(t)\,Jn$ for $m\in\mathbb S^2\setminus\{\pm e_3\}$ 
and $\Phi(\pm e_3)=0$, so that $T=\Phi\circ m$. 
The unit tangent vectors $\tau_1$ and $\tau_2$ of Section~\ref{sec:calibration} correspond to $\partial_t$ and $(\sin t)^{-1}\partial_\phi$. 
Since $\partial_\phi(Jn)=-n$, the derivatives of $\Phi$ along them are
\[
D\Phi(m)\tau_1=F_t\,Jn,\qquad D\Phi(m)\tau_2=-\frac Fs\,n .
\]
Hence $|D\Phi(m)|\le C|m-e_3|$ for $m\ne\pm e_3$; in particular, $D\Phi$ is bounded. 
Since $F$ vanishes at both poles, $\Phi$ is continuous on $\mathbb S^2$, and therefore Lipschitz. 
As $\Phi(e_3)=0$, it follows that $|T|\le C|m-e_3|$, so $T\in L^2$.

\emph{$\diver T\in L^1$.} As for $t$, the chain rule gives $T\in H^1_{\rm loc}$, 
with $\nabla T=D\Phi(m)\nabla m$ almost everywhere on $\Omega$ and $\nabla T=0$ almost everywhere on the pole sets. 
Hence $|\nabla T|\le C|m-e_3||\nabla m|$ almost everywhere, 
and $\nabla T\in L^1$ by the Cauchy--Schwarz inequality. 
Moreover, inserting $\partial_im=u_i\tau_1+v_i\tau_2$ into the chain rule gives \eqref{fluxlocal} almost everywhere on $\Omega$.

\emph{$\int\diver T=0$.} Let $\chi_R$ be the cutoff function in the proof of Lemma~\ref{weighteddegree}. 
Then $\|\nabla\chi_R\|_{L^2}^2\le(C/R)^2|B_{2R}|$ is bounded independently of $R$, so
\[
\Big|\int\chi_R\,\diver T\Big|=\Big|\int T\cdot\nabla\chi_R\Big|
\le\|\nabla\chi_R\|_{L^2}\,\|T\|_{L^2(\{|x|>R\})}\longrightarrow0 .
\]
By the dominated convergence theorem, the left side tends to $\big|\int\diver T\big|$, which is therefore zero.

\emph{The energy identity.} Let $m\in\Aminus$. 
Since \eqref{fluxlocal} holds almost everywhere on $\Omega$ and the proof of Lemma~\ref{lem:calibration} is algebraic, \eqref{CAL} holds almost
everywhere on $\Omega$, with $\det(u,v)=j_m$ and $\diver(F\,Jn)=\diver T$. 
Each term of \eqref{CAL} is integrable: $e_{\rm cal}$ because $m-e_3\in H^1$; 
$H(t)\,j_m$ because $H$ is bounded and $|j_m|\le\frac12|\nabla m|^2$; $\diver T$ by the above; 
and $\mathcal R_\mu(X,Y)$ because $|\mathcal R_\mu(X,Y)|\le|X|^2+|Y|^2$, 
where $X,Y\in L^2$ since $u,v\in L^2$ and $A(t),\beta(t)\le C|m-e_3|$. 
We integrate \eqref{CAL} over $\Omega$. 
Since $j_m=0$ and $\diver T=0$ almost everywhere on the pole sets, Lemma~\ref{weighteddegree} with $G=H$, 
which is continuous by \eqref{calgap}, and $Q(m)=-1$ gives $\int_\Omega H(t)\,j_m=C_\mu$, and the previous step gives $\int_\Omega\diver T=0$. 
Finally, the energy density of $m$ is $e_{\rm cal}$ almost everywhere on $\Omega$, 
and as noted at the beginning of this section, its integral over the pole sets is $2|\{m=-e_3\}|$. This proves \eqref{globalCAL}.
\end{proof}

For the Bloch field, $X=Y=0$ and $\{B_f=-e_3\}=\{0\}$, so $E_\mu(B_f)=C_\mu$. 
Together with \eqref{remainderpositive}, \eqref{globalCAL} gives, for every $m\in\Aminus$,
\begin{equation}
E_\mu(m)-E_\mu(B_f)\ge\frac12\int_\Omega\big(1-|H(t)|\big)\big(|X|^2+|Y|^2\big)\dd x
+2|\{m=-e_3\}|.                                                \label{comparison}
\end{equation}
By \eqref{calgap}, the right side is at least $\frac{\delta_\mu}{2}\int_\Omega(|X|^2+|Y|^2)\dd x+2|\{m=-e_3\}|$.
In particular, $B_f$ minimizes $E_\mu$ on $\Aminus$. 
In the next section, we determine all $m$ for which equality holds.

\section{The equality case}
\label{sec:equality}

\begin{proposition}\label{equality}
Let $m\in\Aminus$. Then $E_\mu(m)=E_\mu(B_f)$ 
if and only if $m=B_f(\cdot-a)$ almost everywhere for some $a\in\R^2$.
\end{proposition}
\begin{proof}
The energy is invariant under translations, so every translate of $B_f$ has energy $E_\mu(B_f)$. 
Conversely, let $E_\mu(m)=E_\mu(B_f)$. 
By \eqref{comparison}, $\int_\Omega(1-|H(t)|)(|X|^2+|Y|^2)\dd x=0$. 
Since $0<t<\pi$ on $\Omega$ and $|H(t)|<1$ for $0<t<\pi$ by \eqref{calgap}, we get $X=Y=0$ almost everywhere on $\Omega$; that is,
\begin{equation}\label{firstorder}
u=A(t)\,Jn,\qquad v=\beta(t)\,n .
\end{equation}

\emph{Regularity.} By \eqref{firstorder}, $|\nabla m|^2=A^2+\beta^2$ almost everywhere on $\Omega$, 
and $A$ and $\beta$ are bounded by Section~\ref{sec:target}; on the pole sets, $\nabla m=0$ almost everywhere. 
Hence, $\nabla m\in L^\infty$ and $m$ has a Lipschitz continuous representative, which we use from now on. 
In particular, $\Omega$ is open. On $\Omega$, \eqref{firstorder} and $\partial_im=u_i\tau_1+v_i\tau_2$ give
\[
\partial_im=A(t)\,(Jn)_i\,\tau_1+\beta(t)\,n_i\,\tau_2\qquad\text{almost everywhere}.
\]
Since $f'<0$, the functions $R=f^{-1}$, $A$ and $\beta$ are smooth on $(0,\pi)$. 
Hence, the right side is a smooth function of $m\in\mathbb S^2\setminus\{\pm e_3\}$, and therefore continuous on $\Omega$. 
A Lipschitz function whose derivatives agree almost everywhere with continuous functions is $C^1$. 
Thus $m\in C^1(\Omega)$, and \eqref{firstorder} holds everywhere on $\Omega$.

\emph{The center.} The set $\Omega$ is nonempty, since otherwise $j_m=0$ almost everywhere and $Q(m)=0$. 
Let $U$ be a connected component of $\Omega$, and write $R$ for $R(t(x))$ on $U$. 
Since $R_t=-1/A$ by Section~\ref{sec:target}, \eqref{firstorder} gives
$\nabla R=R_t\nabla t=-Jn$. Since $|n|=1$, we have $\partial_in=(Jn\cdot\partial_in)\,Jn$; 
moreover, $w=\sin t\,n$ gives $v_i=Jn\cdot\partial_iw=\sin t\,Jn\cdot\partial_in$, and $\beta=\sin t/R$. 
Hence \eqref{firstorder} gives $\partial_in=(n_i/R)\,Jn$. Using $e_i=n_i\,n+(Jn)_i\,Jn$ and $J^2=-1$, we obtain
\[
\partial_i(x+R\,Jn)=e_i-(Jn)_i\,Jn+R\,J\Big(\frac{n_i}{R}\,Jn\Big)=e_i-(Jn)_i\,Jn-n_i\,n=0 .
\]
Since $U$ is connected, $x+R\,Jn=a$ on $U$ for some $a\in\R^2$. As $|Jn|=1$ and $R>0$, this gives
\[
R=|x-a|,\qquad n=J\,\frac{x-a}{|x-a|},\qquad t=f(|x-a|)\qquad\text{on }U .
\]
Thus $a\notin U$. Since $e_\theta=(Jx/|x|,0)$, we also get $m=B_f(\cdot-a)$ on $U$.

\emph{$U=\R^2\setminus\{a\}$.} Let $y\in\partial U$ with $y\ne a$. 
By continuity, $m(y)=B_f(y-a)$, whose polar angle $f(|y-a|)$ lies in $(0,\pi)$; thus $y\in\Omega$. 
Since $\Omega$ is open, a small disk around $y$ lies in $\Omega$. 
This disk is connected and meets $U$, so it lies in $U$, which contradicts $y\in\partial U$ because $U$ is open. 
Hence, $U$ is open and closed in the connected set $\R^2\setminus\{a\}$, so $U=\R^2\setminus\{a\}$. 
Finally, by continuity, $m(a)=\lim_{x\to a}B_f(x-a)=-e_3$. Thus $m=B_f(\cdot-a)$ everywhere.
\end{proof}

We summarize what has been proved. 
Only the hypothesis of Lemma~\ref{triangles} was used after Section~\ref{sec:barriers}.

\begin{theorem}[Conditional classification]\label{conditional}
Let $\mu>0$ and let $f$ be a profile in \eqref{profileode} with $a_0=-f'(0)\ge1$, where $\mu>3/4$ if $a_0=1$. 
Then
\[
\arg\min_{\Aminus}E_\mu=\{B_f(\cdot-a):a\in\R^2\},
\]
and every other profile $\tilde f\ne f$ in \eqref{profileode} satisfies $E_\mu(B_{\tilde f})>E_\mu(B_f)$. 
In particular, $f$ is the only profile in \eqref{profileode} with $a_0\ge1$.
\end{theorem}

\begin{proof}
The first assertion is \eqref{comparison} and Proposition~\ref{equality}. 
If $\tilde f$ is another profile, then $B_{\tilde f}\in\Aminus$; 
so $E_\mu(B_{\tilde f})\ge E_\mu(B_f)$, with equality only if $B_{\tilde f}=B_f(\cdot-a)$ for some $a$. 
The map $B_{\tilde f}$ takes the value $-e_3$ only at $0$, and $B_f(\cdot-a)$ only at $a$; hence $a=0$ and $\tilde f=f$.
Finally, if a profile $\tilde f\ne f$ had $\tilde a_0=-\tilde f'(0)\ge1$, 
then the argument would apply to $\tilde f$ as well (if $\tilde a_0=1$, then $\mu>3/4$, since otherwise $\mu\le1$ and Lemma~\ref{originslope} gives $\tilde a_0>1$), 
and would give $E_\mu(B_f)>E_\mu(B_{\tilde f})$, contradicting what was just proved.
\end{proof}

For $0<\mu\le1$, Proposition~\ref{profileexists} gives a profile, 
and Lemma~\ref{originslope} gives $a_0>1$ for every profile. 
Hence, for $0<\mu\le1$, the profile in \eqref{profileode} is unique, 
and $\arg\min_{\Aminus}E_\mu$ consists of the translates of its Bloch field.

\begin{proposition}[Sharpness of the slope condition]\label{sharpness}
Let $\mu>0$ and let $f$ be a profile in \eqref{profileode} with $a_0=-f'(0)\ge1$. 
Then $\mu<4/\pi$. Moreover, if $\mu>4/\pi$, then $\inf_{\Aminus}E_\mu=-\infty$.
\end{proposition}

\begin{proof}
Let $\theta(x)=4\arctan e^{-x}$. Then $\theta$ decreases from $2\pi$ to $0$, 
and $\theta'=-2\sin(\theta/2)$, so that $\frac12\theta'^2=1-\cos\theta$. 
Consider the wall $q(x)=(0,\sin\theta(x_1),\cos\theta(x_1))$. 
With the DMI density written as in \eqref{dmiext}, its energy per unit length in the $x_2$ direction is
\[
T=\int_{\R}\Big(\frac12\theta'^2+2\mu\sin^2\theta\,\theta'+1-\cos\theta\Big)\dd x_1=8-2\pi\mu .
\]
We compute $T$ in a second way. Since $q$ depends on $x_1$ only, $j_q=0$. 
For $x_1\ne0$ we have $q\ne\pm e_3$, and \eqref{CAL} gives $e_{\rm cal}=\Phi'(x_1)+\mathcal R_\mu(X,Y)$, 
where $\Phi=F(t)\,(Jn)_1$. 
The function $\Phi$ is continuous on $\R$ because $F(\pi)=0$, 
and it vanishes at $\pm\infty$ because $F(0)=0$, and $\Phi\in W^{1,1}(\R)$, since $|\Phi'|\le\sup|F_t|\,|\theta'|=2\sup|F_t|\,\sin(\theta/2)$ is integrable.
Hence $\int_\R\Phi'=0$ and $T=\int_\R\mathcal R_\mu(X,Y)\dd x_1$. 
Along the wall, $n$ is constant on each half-line, so $v=0$ and $Y=-\beta(t)\,n$, where $\beta(t)>0$ for $0<t<\pi$. 
If $\mu\le3/4$, there is nothing to prove, so let $\mu>3/4$. 
Since $a_0\ge1$, Lemma~\ref{triangles} gives $|H(t)|<1$ for $0<t<\pi$, 
so $\mathcal R_\mu(X,Y)\ge\frac12(1-|H(t)|)\beta(t)^2>0$ for $x_1\ne0$. Thus $8-2\pi\mu=T>0$.

For the second claim, let $g=\theta|_{[0,\infty)}$, which is a decreasing function from $\pi$ to $0$, 
and let $B_g$ be the Bloch field with profile $g$. 
On the line $x_2=0$, $B_g$ coincides with $q$. 
For $L>0$, define $m_L=q$ on $\{|x_2|\le L/2\}$ and $m_L(x)=B_g(x\mp(0,L/2))$ on $\{\pm x_2>L/2\}$. 
Then $m_L$ is continuous and $m_L-e_3\in H^1$, 
the two half-planes together carry the energy and the degree of $B_g$, and the strip carries no degree. 
Hence $m_L\in\Aminus$ and $E_\mu(m_L)=E_\mu(B_g)+(8-2\pi\mu)L$, which tends to $-\infty$ as $L\to\infty$ when $\mu>4/\pi$.
\end{proof}

\section{The radial branch and the calibration range}
\label{sec:branch}

In this section we use Appendix~\ref{ex-existence},
which provides, for every $0<\mu\le\sqrt2$, the existence of profiles, 
the compactness of the set of profiles over compact parameter intervals (Proposition~\ref{ex-proper}), 
and the fact that near every profile the set of profiles is the graph of a smooth function of $\mu$ (the proof of Proposition~\ref{ex-open}).

\begin{proposition}[The radial branch]\label{branch}
For every $0<\mu\le\sqrt2$, the equation \eqref{profileode} has exactly one solution $f_\mu$. 
The map $\mu\mapsto f_\mu-f_{\mu_0}$ is smooth near every $\mu_0$ with values in the space $\mathcal H_r$ of Lemma~\ref{ex-scalar}, 
and $a(\mu)=-f_\mu'(0)$ is a smooth function on $(0,\sqrt2]$.
\end{proposition}

\begin{proof}
Let $N(\mu)$ be the number of profiles for the parameter $\mu$. 
Each profile is an isolated point of its fiber, since near it the profiles form a graph over $\mu$, 
and the fiber is compact by Proposition~\ref{ex-proper}; hence $N(\mu)<\infty$. 
We claim that $N$ is locally constant on $(0,\sqrt2]$. 
Fix $\mu_0$, let $f^1,\dots,f^N$ be the profiles at $\mu_0$, 
and choose disjoint neighborhoods $V_1,\dots,V_N$ in which the profiles form graphs over $\mu$. 
For $\mu$ close to $\mu_0$, each $V_i$ contains exactly one profile, 
so $N(\mu)\ge N$. If $N(\mu_k)>N$ for some $\mu_k\to\mu_0$, there are profiles at $\mu_k$ outside $\bigcup V_i$; 
by Proposition~\ref{ex-proper}, 
a subsequence converges to a profile at $\mu_0$ outside $\bigcup V_i$, which is impossible. 
Since $(0,\sqrt2]$ is connected and $N=1$ on $(0,1]$ by Proposition~\ref{profileexists}, Lemma~\ref{originslope}, and Theorem~\ref{conditional}, we get $N\equiv1$.
The smoothness statements follow from the graph structure and from the smooth dependence of the origin shots 
and of the tails on their data (Lemmas~\ref{ex-origin-lemma} and~\ref{ex-tailgraph}).
\end{proof}

\begin{lemma}[The slope decreases]\label{slopedecreasing}
The function $a(\mu)=-f_\mu'(0)$ is strictly decreasing on $(0,\sqrt2]$.
\end{lemma}

\begin{proof}
Fix $\mu$ and put $w=\partial_\mu f_\mu$. 
By Proposition~\ref{branch}, $w\in\mathcal H_r$; 
near the origin, $f_\mu(0)=\pi$ for all $\mu$ gives $w=O(r)$, 
and at infinity $w$ decays exponentially by Lemma~\ref{ex-tailgraph}. 
Differentiating the equation in \eqref{profileode} with respect to $\mu$ gives
\[
L_fw=\frac{2\sin^2f}{r}=:g>0,\qquad
L_f=-\partial_r^2-r^{-1}\partial_r+V_f,\qquad V_f=\frac{\cos2f}{r^2}-\frac{2\mu\sin2f}{r}+\cos f,
\]
and the equation holds weakly against every test function in $\mathcal H_r$, 
since $g\in L^2(r\dd r)$ and $|V_f|\le C(1+r^{-2})$. 
Let $w_-=\max(-w,0)\in\mathcal H_r$. 
Testing with $w_-$ gives $-q_f[w_-]=\int_0^\infty g\,w_-\,r\dd r\ge0$, where $q_f[v]=\int_0^\infty(v'^2+V_fv^2)\,r\dd r$. 
By Lemma~\ref{ex-scalar}, $q_f[v]>0$ for $v\ne0$; hence $w_-=0$ and $w\ge0$. 
If $w(r_0)=0$ for some $r_0>0$, then $w'(r_0)=0$ and $w''(r_0)\ge0$, so $L_fw(r_0)\le0<g(r_0)$, a contradiction. 
Thus $w>0$ on $(0,\infty)$. 
Differentiating \eqref{originseries} with respect to $\mu$, which is justified by the smooth dependence in Lemma~\ref{ex-origin-lemma}, gives
\[
w(r)=-a'(\mu)\,r+\Big(\partial_\mu c_3+\partial_{a_0}c_3\,a'(\mu)\Big)r^3+O(r^5),
\qquad \partial_\mu c_3=-\frac{a_0^2}{4}.
\]
Since $w>0$, we get $a'(\mu)\le0$; and $a'(\mu)=0$ would give $w=-\frac{a_0^2}4r^3+O(r^5)<0$ near the origin. Hence $a'(\mu)<0$.
\end{proof}

\begin{proposition}[The calibration endpoint]\label{mucal}
There is exactly one $\mu_{\rm cal}\in(1,4/\pi)$ with $a(\mu_{\rm cal})=1$. 
One has $a(\mu)\ge1$ for $0<\mu\le\mu_{\rm cal}$ and $a(\mu)<1$ for $\mu_{\rm cal}<\mu\le\sqrt2$. 
For $0<\mu\le\mu_{\rm cal}$,
\[
\arg\min_{\Aminus}E_\mu=\{B_{f_\mu}(\cdot-a):a\in\R^2\}.
\]
\end{proposition}

\begin{proof}
By Lemma~\ref{originslope}, $a(1)>1$. 
By Proposition~\ref{sharpness}, a profile with $a_0\ge1$ exists only for $\mu<4/\pi$; 
since $4/\pi<\sqrt2$, this gives $a(4/\pi)<1$. As $a$ is continuous and strictly decreasing (Lemma~\ref{slopedecreasing}), 
there is exactly one $\mu_{\rm cal}\in(1,4/\pi)$ with $a(\mu_{\rm cal})=1$, and the two inequalities follow. 
For $0<\mu\le\mu_{\rm cal}$ we have $a(\mu)\ge1$ and $\mu>3/4$ whenever $a(\mu)=1$, so Theorem~\ref{conditional} applies to $f_\mu$.
\end{proof}

\begin{remark}
Numerically, $\mu_{\rm cal}\approx1.1714$, that is, $\mu_{\rm cal}^{-2}\approx0.729$; 
the energy $E_\mu(B_{f_\mu})$ vanishes at $\mu\approx1.117$, that is, at $h/\kappa^2\approx0.80$, 
which is the field below which the skyrmion lattice has lower energy than the uniform state~\cite{BH94}. 
These values are numerical observations and are not used in the proofs.
\end{remark}

\section{Proofs of the main theorems}
\label{sec:original}

We now prove Theorems~\ref{originaltheorem} and~\ref{rangetheorem} by reducing \eqref{originalenergy} to the normalized energy $E_\mu=E_{\mu,1}$, 
which we use here for every $\mu\in\R$.

Let $\kappa\ne0$ and $h>0$. For $m\in\Aminus$, put $M(y)=m(y/\sqrt h)$, that is, $m(x)=M(\sqrt h\,x)$. 
The dilation $y=\sqrt h\,x$ preserves orientation, 
and a change of variables using the last integral in \eqref{dmiext} for $D$ gives
\[
\int|\nabla m|^2=\int|\nabla M|^2,\quad D(m)=h^{-1/2}D(M),\quad
h\int|m-e_3|^2=\int|M-e_3|^2,\quad Q(M)=Q(m).
\]
Hence $m\mapsto M$ is a bijection of $\Aminus$ onto itself, and
\begin{equation}\label{scaling}
E_{\kappa,h}(m)=E_{\kappa/\sqrt h}(M).
\end{equation}

Let $S=\operatorname{diag}(-1,-1,1)$ be the rotation by $\pi$ about $e_3$. 
The exchange and Zeeman energies are invariant under $m\mapsto Sm$, and so is the degree, since $S\in SO(3)$; 
by \eqref{dmiext}, $D(Sm)=-D(m)$. 
Hence $m\mapsto Sm$ is a bijection of $\Aminus$ onto itself, and
\begin{equation}\label{signchange}
E_{-\nu}(Sm)=E_\nu(m)\qquad(\nu\in\R).
\end{equation}

\begin{proof}[Proof of Theorem~\ref{originaltheorem}]
Let $h\geq \kappa^2$, so that $\nu=|\kappa|/\sqrt h\in(0,1]$. 
By Proposition~\ref{profileexists} and the uniqueness argument following Theorem~\ref{conditional}, 
there is exactly one profile $f_\nu$ in \eqref{profileode} with $\mu=\nu$, and $\arg\min_{\Aminus}E_\nu=\{B_{f_\nu}(\cdot-b):b\in\R^2\}$.

If $\kappa>0$, then by \eqref{scaling}, 
$m$ minimizes $E_{\kappa,h}$ on $\Aminus$ if and only if $M=B_{f_\nu}(\cdot-b)$ almost everywhere for some $b\in\R^2$. 
If $\kappa<0$, then $E_{\kappa/\sqrt h}(M)=E_{-\nu}(M)=E_\nu(SM)$ by \eqref{signchange} and $S^2=I$, 
so $m$ is a minimizer if and only if $SM=B_{f_\nu}(\cdot-b)$ almost everywhere for some $b$. 
Since $B_f(x)=\big(\sin f(|x|)\,Jx/|x|,\cos f(|x|)\big)$, 
in both cases $m$ is a minimizer if and only if $m=m_{\kappa,h,a}$ almost everywhere for some $a\in\R^2$ (namely $a=b/\sqrt h$), 
with $m_{\kappa,h,a}$ as in \eqref{originalfamily}.

Each $m_{\kappa,h,a}$ is smooth. Away from $a$ this holds because $f_\nu$ is smooth on $(0,\infty)$. 
Near $a$, it holds because $\sin f_\nu(r)/r$ and $\cos f_\nu(r)$ are smooth functions of $r^2$ by the proof of Lemma~\ref{ex-origin-lemma}; 
thus the planar component of $m_{\kappa,h,a}$ is a smooth function of $|x-a|^2$ times $J(x-a)$, 
and the vertical component is a smooth function of $|x-a|^2$.

Since $0<f_\nu<\pi$ on $(0,\infty)$, the map $m_{\kappa,h,a}$ takes the value $-e_3$ only at $a$. 
As these maps are continuous, two of them agree almost everywhere only if they have the same $a$. 
Hence every minimizer agrees almost everywhere with exactly one $m_{\kappa,h,a}$. 
Finally, let $m=m_{\kappa,h,a}$. 
Since $|R_\theta y|=|y|$ and $JR_\theta=R_\theta J$, \eqref{originalfamily} gives $m(a+R_\theta y)=\widehat R_\theta\,m(a+y)$. 
Conversely, if $m(c+R_\theta y)=\widehat R_\theta\,m(c+y)$ for all $\theta$ and $y$, 
then taking $y=a-c$ gives $m(c+R_\theta(a-c))=\widehat R_\theta(-e_3)=-e_3$. 
Hence $c+R_\theta(a-c)=a$ for all $\theta$, and $c=a$. Thus $a$ is the only center of this rotational symmetry.
\end{proof}

\begin{proof}[Proof of Theorem~\ref{rangetheorem}]
The assertions about the profiles are Proposition~\ref{branch}, Lemma~\ref{slopedecreasing} and Proposition~\ref{mucal}. 
Let $\kappa\ne0$ and $h>0$ with $\nu=|\kappa|/\sqrt h\le\mu_{\rm cal}$. 
By Proposition~\ref{mucal}, $f_\nu$ is the only profile in \eqref{profileode} with $\mu=\nu$, and $\arg\min_{\Aminus}E_\nu=\{B_{f_\nu}(\cdot-b):b\in\R^2\}$; 
the proof of Theorem~\ref{originaltheorem} now applies without change. Finally, let $h<\pi^2\kappa^2/16$, that is, $\nu>4/\pi$. 
By \eqref{scaling} and \eqref{signchange}, $\inf_{\Aminus}E_{\kappa,h}=\inf_{\Aminus}E_\nu$, which is $-\infty$ by Proposition~\ref{sharpness}.
\end{proof}

\appendix

\providecommand{\Pmu}{\mathcal P_\mu}
\providecommand{\Hr}{\mathcal H_r}
\section{Existence of a monotone Bloch profile on the whole interval
  \texorpdfstring{$0<\mu \leq \sqrt{2}$}{0 < mu <= sqrt(2)}}\label{ex-existence}

\noindent
This appendix proves Proposition~\ref{profileexists}: for every $0<\mu \leq \sqrt{2}$ there is a profile satisfying \eqref{profileode}. 
The proof is a continuation argument in $\mu$. 
Sections~A.1 and~A.2 derive asymptotics and a priori estimates valid for every such profile. 
Using them, we show that the set of $\mu\in(0,\sqrt{2}]$ for which such a profile exists is closed (Section~A.3) and open (Section~A.4) in $(0,\sqrt{2}]$. 
By a result of Ibrahim--Shimizu~\cite{IS}, it contains $\mu=1/2$ (Section~A.5), and hence it equals $(0,\sqrt{2}]$. 
The argument does not use the calibration, minimality, or symmetry results of the main text.
The notation of this appendix is self-contained.

\subsection{Main result and asymptotics}

The curl convention is
\[
 \operatorname{curl}m=(\partial_2m_3,-\partial_1m_3,
                         \partial_1m_2-\partial_2m_1).
\]
The energy is
\[
 E_\mu(m)=\frac12\int_{\R^2}|\nabla m|^2
 +\mu\int_{\R^2}(m-e_3)\cdot\operatorname{curl}(m-e_3)
 +\frac12\int_{\R^2}|m-e_3|^2.
\]
For $m_f=\sin f\,e_\theta+\cos f\,e_3$, its radial Euler equation is
\begin{equation}
 -f''-\frac{f'}r+\frac{\sin f\cos f}{r^2}
       -\frac{2\mu\sin^2f}{r}+\sin f=0.               \label{ex-ode}
\end{equation}
Indeed integration by parts gives the radial energy divided by $2\pi$ as
\[
 E_e(f)+\mu D(f)+Z(f),\quad
 E_e(f)=\frac12\int_0^\infty
             \left(f'^2+\frac{\sin^2f}{r^2}\right)r\,dr,
\]
\[
 D(f)=2\int_0^\infty r f'\,\sin^2f \,dr,\qquad
 Z(f)=\int_0^\infty(1-\cos f)r\,dr.
\]
The DMI expression follows equally from
$D(m)=-2\int m_\parallel\cdot J\nabla m_3$, where
$J(x_1,x_2)=(-x_2,x_1)$, so it uses the indicated helicity and curl sign.

Let $\Pmu$ consist of all smooth solutions of \eqref{ex-ode} on
$(0,\infty)$ such that
\begin{equation}
 f(0)=\pi,\quad f(\infty)=0,\quad 0<f<\pi,\quad f'<0,
 \qquad E_e(f)+Z(f)<\infty.                         \label{ex-class}
\end{equation}

\begin{theorem}\label{ex-main}
For every fixed $0<\mu \leq \sqrt{2}$ the class $\Pmu$ is nonempty. Each member
has a positive finite origin slope and the differentiated origin
expansion \eqref{ex-origin}. At infinity it has the refined positive
massive asymptotics \eqref{ex-tail}.
\end{theorem}

\begin{lemma}[Origin asymptotics and origin shots]\label{ex-origin-lemma}
Every $f\in\Pmu$ has $a=-f'(0)\in(0,\infty)$ and
\begin{equation}
 f(r)=\pi-ar+
       c_3r^3+O(r^5), \quad c_3:=\frac{a+(2/3)a^3-2\mu a^2}{8}.     \label{ex-origin}
\end{equation}
The differentiated remainder is $O(r^4)$. Conversely, the data $(\mu,a)$
determine a unique regular origin solution, with smooth dependence and
uniform local remainders on compact parameter and slope sets.
\end{lemma}
\begin{proof}
Set $t=\log r$ and $g=\pi-f$. Then $g>0$, $g_t=-r f^{\prime}>0$, $g(-\infty)=0$, and
\[
 g_{tt}=\sin g\cos g+2\mu e^t\sin^2g-e^{2t}\sin g
        =g(1+h), \quad
 h=O(g^2+e^tg+e^{2t})\longrightarrow0 \quad \text{ as } t \rightarrow -\infty.
\]
Fix $T$ with $|h| \leq 1 / 4$ on $(-\infty, T]$; all estimates below are taken there. Then
$3g/4\le g_{tt}\le5g/4$, so $g$ is convex and $g_t$ is increasing. Convexity and $g(-\infty)=0$ imply
$g_t(-\infty)=0$. Multiplying by $g_t$ and integrating gives
$\sqrt3/2\le q=g_t/g\le\sqrt5/2$. Hence
$g=O(e^{\alpha t})$, $\alpha=\sqrt3/2$, and
$h=O(e^{2\alpha t})$. 

Differentiation of $q$ gives the Riccati equation $q_t=g_{t t} / g-q^2=1-q^2+h$. Put $d:=q-1$, 
so $|d| \leq 1-\alpha$ by the two-sided bound on $q$. 
Using the identity $1-q^2=-(q+1) d$,
the bounded function $d$ satisfies
\[
 d_t+(q+1)d=h,\qquad q+1 \geq 1+\alpha>0.
\]
For $s<t \leq T$, variation of constants gives the identity
\[
d(t)=\exp \left(-\int_s^t(q+1)\right) d(s)+\int_s^t \exp \left(-\int_\sigma^t(q+1)\right) h(\sigma) d \sigma .
\]
The kernel is bounded by $e^{-(1+\alpha)(t-\sigma)}$, so the first term tends to 0 as $s \rightarrow-\infty$ since $d$ is bounded, and, as $|h(\sigma)| \leq C e^{2 \alpha \sigma}$ with $1+3 \alpha>0$, the second term converges absolutely. Therefore
\[
d(t)=\int_{-\infty}^{t}\exp\Bigl(-\int_{s}^{t}(q+1)\Bigr)h(s)ds,\qquad
|d(t)|\le C e^{-(1+\alpha)t}\int_{-\infty}^{t}e^{(1+3\alpha)s}ds=\frac{C}{1+3\alpha}e^{2\alpha t} .
\]
This is the only bounded solution: the omitted homogeneous solution $\exp (-\int_T^t(q+1 ))$ grows at least as fast as $e^{(1+\alpha)|t|}$ at minus infinity.
Since $d=\partial_t\log(e^{-t}g)$ is integrable there, $\log \left(e^{-t} g\right)$ converges and
$e^{-t}g\to a\in(0,\infty)$.

Write $f=\pi-ru$. Its equation becomes
\[
 u''+3u'/r=F(r,u,\mu),\quad
 F=\frac{\sin(ru)\cos(ru)-ru}{r^3}
       +2\mu\frac{\sin^2(ru)}{r^2}-\frac{\sin(ru)}r.
\]
All singularities of $F$ are removable: it is analytic and even in $r$,
with $F(0,u,\mu)=-u-(2/3)u^3+2\mu u^2$. 
The left side is in divergence form, $u^{\prime \prime}+3 u^{\prime} / r=r^{-3}\left(r^3 u^{\prime}\right)^{\prime}$, and its flux vanishes at the origin. Indeed $u=g / r$ gives, with $q=g_t / g=r g_r / g$,
\[
r^3 u^{\prime}=r^2 g_r-r g=r g(q-1),
\]
where $|q-1| \leq 1-\alpha$ by the preceding two-sided bound and $r g \rightarrow 0$; the sharper estimates $g=O(r)$ and $q-1=O\left(r^{2 \alpha}\right)$ give $r^3 u^{\prime}=O\left(r^{2+2 \alpha}\right)$ and $u^{\prime}=O\left(r^{2 \alpha-1}\right)$ with $2 \alpha-1>0$. Hence integration of $\left(r^3 u^{\prime}\right)^{\prime}=r^3 F$ over $(0, r)$ leaves no boundary term,
\[
u^{\prime}(r)=r^{-3} \int_0^r \sigma^3 F(\sigma, u(\sigma), \mu) d \sigma
\]
and since $u=g / r \rightarrow a$ while $u^{\prime}$ is integrable at the origin, a second integration over $(0, r)$ gives
\begin{equation}
 u(r)=a+\int_0^r t^{-3}\int_0^t
                         \tau^3F(\tau,u(\tau),\mu)\,d\tau\,dt.
                                                        \label{ex-volterra}
\end{equation}
Work in $C\left(\left[0, r_0\right]\right)$, on the unit ball about $a$, where $|u| \leq a+1$; 
let $M$ bound $|F|$ and $L$ bound $\left|F_u\right|$ there. Since $\int_0^r t^{-3} \int_0^t \tau^3 d \tau d t=r^2 / 8$, 
the right side of \eqref{ex-volterra} maps that ball into itself if $M r_0^2 / 8 \leq 1$ and is contractive if $L r_0^2 / 8<1$; 
choose $r_0$ accordingly, so that \eqref{ex-volterra} has a unique solution there. 
As $M$ and $L$ depend on the data only through $a$ and the range of $\mu$, a single $r_0$ serves each compact data set.

To establish regularity at every
order, set $p=r^2$, $U(p)=u(\sqrt p)$, and
$G(p,U)=F(\sqrt p,U,\mu)$. Since $F$ is even analytic in $r$, $G$ is
smooth at $p=0$. Differentiating \eqref{ex-volterra} for $p>0$ gives
\[
 U'(p)=\frac14\int_0^1 t\,G(pt,U(pt))\,dt.
\]
The right side is continuous up to zero when $U$ is continuous, so $U$
is $C^1$ up to zero. Differentiating this formula then bootstraps $U$
to $C^\infty$. Thus $u(r)$ is a smooth function of $r^2$. 
In particular,
$\sin f(r)/r=U(r^2)\operatorname{sinc}(rU(r^2))$ and
$\cos f(r)=-\cos(rU(r^2))$ are smooth functions of $r^2$;
the Cartesian Bloch field is smooth through its center. Here
$\operatorname{sinc}x=\sin x/x$ is continued smoothly at $x=0$.
Differentiation of the fixed-point equation and its parameter equations
also gives smooth data dependence. 
Since $u(r)=U\left(r^2\right)$ with $U$ smooth at 0 and $U(0)=a$, the last display at $p=0$ gives
\[
U^{\prime}(0)=\frac{1}{4} G(0, a) \int_0^1 \theta d \theta=\frac{1}{8} F(0, a, \mu)=-c_3, \quad c_3:=\frac{a+(2 / 3) a^3-2 \mu a^2}{8}.
\]
Taylor expansion at $p=0$ gives $U(p)=a-c_3 p+O\left(p^2\right)$ and $U^{\prime}(p)=-c_3+O(p)$, hence, with $p=r^2$,
\[
u=a-c_3 r^2+O\left(r^4\right), \quad u^{\prime}=2 r U^{\prime}\left(r^2\right)=-2 c_3 r+O\left(r^3\right) .
\]
Since $f=\pi-r u$ and $f^{\prime}=-u-r u^{\prime}$, this is \eqref{ex-origin} together with $f^{\prime}=-a+3 c_3 r^2+O\left(r^4\right)$.
The original profile satisfies this same fixed-point equation and, 
after shrinking $r_0$, lies in that ball; it is therefore the fixed point, proving uniqueness at the origin. On any interval 
bounded away from zero the trigonometric forcing is bounded, so $\left|f^{\prime \prime}\right| \leq C\left(1+\left|f^{\prime}\right|\right)$ 
and each origin shot has ordinary finite-radius continuation and smooth dependence on its data.
\end{proof}

\begin{lemma}[Asymptotics at infinity]\label{ex-tail-lemma}
For every $f\in\Pmu$ there is $c_f>0$ such that, writing $K_1$ for the
modified Bessel function,
\begin{equation}
 f=c_fK_1(r)+o(K_1(r)/r^2),\qquad
 f'=c_fK_1'(r)+o(K_1(r)/r^2).                    \label{ex-tail}
\end{equation}
\end{lemma}
\begin{proof}
For $\eta>0$ set
$B_\eta=-\partial_r^2-r^{-1}\partial_r+r^{-2}+\eta^2$.
Exactly,
\begin{equation}
 B_1f=N_\mu(r,f),\quad
 N_\mu(r,u)=\frac{u-\sin u\cos u}{r^2}
                 +\frac{2\mu\sin^2u}{r}+u-\sin u.  \label{ex-source}
\end{equation}
Fix $1/2<\eta<1$ and choose $R\ge1$ so large that
$N_\mu(r,f)/f<1-\eta^2$ on $r\ge R$. This is possible because
$f$ decreases to zero. Then $B_\eta f<0$. 
Here $I_1$ and $K_1$ denote the modified Bessel functions of order one, 
which are the solutions of $x^2 y^{\prime \prime}+x y^{\prime}-\left(x^2+1\right) y=0$ that are regular at $x=0$ and decay at infinity, respectively. 
Both are positive on $(0, \infty)$, with $I_1(x) \sim e^x / \sqrt{2 \pi x}$ and $K_1(x) \sim \sqrt{\pi /(2 x)} e^{-x}$ as $x \rightarrow \infty$. 
By scaling, $I_1(\eta r)$ and $K_1(\eta r)$ solve $B_\eta y=0$.

The one-dimensional maximum principle gives $f(r) \leq f(R) K_1(\eta r) / K_1(\eta R)$ for $r \geq R$. 
Indeed, fix $\varepsilon>0$ and let $w=f-f(R) K_1(\eta r) / K_1(\eta R)-\varepsilon I_1(\eta r)$. 
Then $B_\eta w=B_\eta f<0$ for $r>R$. Moreover, $w(R)<0$, and $w(\rho)<0$ for all large $\rho$ because $f<\pi$ while $I_1$ grows exponentially. 
For such $\rho$ we have $w \leq 0$ on $[R, \rho]$; 
otherwise, $w$ would attain a positive maximum at some $r_* \in(R, \rho)$, where $w^{\prime}\left(r_*\right)=0$ and $w^{\prime \prime}\left(r_*\right) \leq 0$ so that $B_\eta w\left(r_*\right) \geq\left(r_*^{-2}+\eta^2\right) w\left(r_*\right)>0$, contradicting $B_\eta w<0$. 
Letting $\rho \rightarrow \infty$ and then $\varepsilon \downarrow 0$ gives the claim.

On a sufficiently large tail, writing the right side of $(rf')'=rH$ as
\[
 H=\frac{\sin f\cos f}{r^2}
                   -\frac{2\mu\sin^2f}{r}+\sin f,
\]
we have $0<H\le2f$. 
Thus $r f^{\prime}$ is negative and increasing. 
Its limit is zero, since a negative limit would force logarithmic crossing of zero. 
Integrating $\left(r f^{\prime}\right)^{\prime}=r H$ from infinity gives for $f^{\prime}$ the same submass exponential bound $O\left(K_1(\eta r)\right)$ as for $f$:
\[
0<-r f^{\prime}(r)=\int_r^{\infty} t H(t) d t \leq 2 \int_r^{\infty} t f(t) d t=O\left(r K_1(\eta r)\right).
\]
Here, the mass is the decay rate 1 of the decaying solution $K_1(r) \sim \sqrt{\pi /(2 r)} e^{-r}$ of the linearized equation $B_1 y=0$; 
a bound $O\left(K_1(\eta r)\right)=O\left(r^{-1 / 2} e^{-\eta r}\right)$ with $\eta<1$ decays more slowly and is therefore called submass. 
Finally, $N_\mu$ has no linear term, so $N_\mu(r, f)=O\left(f^2 / r+f^3\right)=O\left(e^{-2 \eta r}\right)$. 
As $r I_1(r)=O\left(r^{1 / 2} e^r\right)$ and $2 \eta>1$, it follows that $\int_R^{\infty} r I_1(r) N_\mu(r, f(r)) d r<\infty$; 
that is, $N_\mu(r, f)$ is integrable against $r I_1(r)$.

Let $I=I_1, K=K_1$ and $H_R(r)=I(r) K(R)-K(r) I(R)$. 
Thus, $K$ is the solution of $B_1 y=0$ that decays at infinity, and $H_R$ is the one that vanishes at $r=R$; 
moreover, $H_R>0$ for $r>R$, since $(I / K)^{\prime}=1 /\left(r K^2\right)>0$. 
For continuous $g$ with $g(r)=O(e^{-\delta r})$ for some $\delta>0$, set
\begin{equation}
 (T_Rg)(r)=\int_R^\infty G_R(r,t)t g(t)\,dt,
 \quad G_R(r,t)=\frac{H_R(\min(r,t))K(\max(r,t))}{K(R)}.
                                                        \label{ex-green}
\end{equation}
For $r \neq t$, $G_R(\cdot, t)$ solves $B_1 y=0$. It is continuous at $r=t$, and 
by the Wronskian identity $I K^{\prime}-I^{\prime} K=-1 / r$, its $r$-derivative jumps 
there by $-1 / t$; with the weight $t$ in \eqref{ex-green}, this jump gives $B_1 T_R g=g$. 
Hence $u=T_R g$ is the solution of $B_1 u=g$ on $(R, \infty)$ with $u(R)=0$ 
and $u \rightarrow 0$ at infinity; this is the Dirichlet decaying Green operator.

Now $v=f-f(R) K / K(R)-T_R N_\mu(f)$ satisfies $B_1 v=0$ and $v(R)=0$, so $v=\beta H_R$ for 
some constant $\beta$. Every term of $v$ 
tends to zero (for $T_R N_\mu(f)$ because $N_\mu(r, f)=O(e^{-2 \eta r})$), whereas 
$H_R$ grows like $I$; 
hence $\beta=0$, that is, $f=f(R) K / K(R)+T_R N_\mu(f)$. 
To isolate the tail, write $\int_R^r=\int_R^{\infty}-\int_r^{\infty}$ in the part of \eqref{ex-green} with $t<r$. 
Since $H_R(r) K(t)-K(r) H_R(t)=$ $K(R)(I(r) K(t)-K(r) I(t))$, this gives
\[
\begin{gathered}
f(r)=c_f K(r)+\int_r^{\infty}(I(r) K(t)-K(r) I(t)) t N_\mu(t, f(t)) d t, \\
c_f=\frac{f(R)}{K(R)}+\frac{1}{K(R)} \int_R^{\infty} H_R(t) t N_\mu(t, f(t)) d t.
\end{gathered}
\]
Thus $c_f$ is the amplitude of the decaying mode $K$ in $f$: the remaining integral sees $N_\mu$ only 
beyond $r$ and is shown below to be $o\left(K(r) / r^2\right)$, so $c_f$ is the constant in \eqref{ex-tail}. 
It is finite by the integrability above, since $0<H_R \leq K(R) I$, and positive because $f(R)>0$, $H_R>0$ and $N_\mu(r, f) \geq 0$.
Denote the difference between $f(r)$ and $c_fK(r)$ by $\mathcal{E}(r)$. For $t \geq r$ its 
kernel satisfies $0 \leq K(r) I(t)-I(r) K(t) \leq K(r) I(t)$, because $I / K$ is 
increasing. With $N_\mu(t, f(t))=O\left(e^{-2 \eta t}\right)$ from 
the submass estimate, $t I(t)=O\left(t^{1 / 2} e^t\right)$ and $K(r)=O\left(r^{-1 / 2} e^{-r}\right)$, this gives
\[
|\mathcal{E}(r)| \leq K(r) \int_r^{\infty} I(t) t N_\mu(t, f(t)) d t=O\left(r^{-1 / 2} e^{-r} \cdot r^{1 / 2} e^{-(2 \eta-1) r}\right)=O\left(e^{-2 \eta r}\right).
\]
Differentiating $\mathcal{E}$ produces no boundary term, because 
the kernel vanishes at $t=r$; since $I^{\prime}$ and $K^{\prime}$ obey 
the same bounds as $I$ and $K$, also $\mathcal{E}^{\prime}(r)=O\left(e^{-2 \eta r}\right)$. 
As $2 \eta>1$, $e^{-2 \eta r}=o\left(K(r) / r^2\right)$, so $f=c_f K+\mathcal{E}$ gives \eqref{ex-tail}.
\end{proof}

\subsection{Stress estimates for every admissible profile}

Lemma \ref{ex-stress} presents a priori bounds for every admissible profile, 
which are used in the next section to prove properness. 
It compares $f$ with the Belavin--Polyakov profiles $2 \arctan (\lambda / r)$, 
which satisfy $-r f^{\prime}=\sin f$ and minimize the exchange energy among degree-one profiles \cite{BP75}. 
The key quantity is the stress $J=r^2 f^{\prime 2}-\sin ^2 f$, which vanishes on Belavin--Polyakov profiles. 
The same quantity, denoted $F$, appears in \cite[Section~3.2]{IS}, 
whose models include the present one ($\beta=1$); 
there, its positivity is proved for the solutions constructed in \cite{IS}, which satisfy the additional constraint \cite[(1.15)]{IS}. 
The argument below is self-contained and applies to every admissible profile.
Using the level $\sigma=\pi-f$ as the independent variable, we show that $J$ increases near the center, 
where the DMI term dominates, and decreases farther out, where the Zeeman term dominates, 
with exactly one switch, the stress transition. Since $J$ vanishes at both ends, 
this gives the strict stress inequality $-r f^{\prime}>\sin f$ and the bound $a>1 /(2 \mu)$ on the origin slope; 
the stress inequality also places $f$ below the Belavin--Polyakov profile with the same origin slope, see \eqref{ex-BP-upper}.
Integrating the equation for $J$ yields $J\le\frac83\mu^2$ 
and the virial identity $\mu D+2Z=0$, from which the energy bounds follow.

\begin{lemma}\label{ex-stress}
For every $0<\mu \leq \sqrt{2}$ and every $f\in\Pmu$, put $A=-rf'$. Then
\begin{equation}
 A>\sin f,\qquad a>\frac1{2\mu},\qquad
 0<A^2-\sin^2f\le\frac83\mu^2,                   \label{ex-stress-bounds}
\end{equation}
and
\begin{equation}
 E_e(f)-2\le\pi\sqrt{2/3}\,\mu,
 \qquad Z(f)\le4\mu^2E_e(f)<48,                  \label{ex-energy-bounds}
\end{equation}
where $a=-f'(0)$.
\end{lemma}
\begin{proof}
Since $f^{\prime}<0$, the level $\sigma=\pi-f(r)$ increases from 0 to $\pi$ as $r$ runs over ($0, \infty$), with $\sigma=0$ 
at the center and $\sigma \rightarrow \pi$ at infinity. We use $\sigma$ 
as the independent variable, write $r=R(\sigma)$ for the inverse function, and 
denote $d / d \sigma$ by a subscript $\sigma$; 
primes still denote $d / d r$. Put $s=\sin \sigma$($=\sin f$) and
\[
J=A^2-s^2=r^2 f^{\prime 2}-\sin ^2 f.
\]
We call $J$ the stress.
It measures the deviation from the Belavin--Polyakov profiles $2 \arctan (\lambda / r)$, which 
satisfy $r f^{\prime}=-\sin f$ and hence $J=0$. 
By \eqref{ex-ode}, $A^{\prime}=-\left(r f^{\prime}\right)^{\prime}=-\sin f \cos f / r+2 \mu \sin ^2 f-r \sin f$, and with $f^{\prime}=-A / r$,
\[
J^{\prime}=2 A A^{\prime}-2 (\sin f \cos f )f^{\prime}=2 A \sin f(2 \mu \sin f-r) .
\]
The exchange terms $\pm 2 A \sin f \cos f / r$ cancel, so $J$ is changed only by the DMI and Zeeman terms; 
without them $J$ would be constant. Since $d \sigma / d r=-f^{\prime}=A / r$, in the variable $\sigma$ this reads
\begin{equation}
R_\sigma=\frac{R}{A}, \quad J_\sigma=2 R s(2 \mu s-R)=2 s\left(\mu^2 s^2-(R-\mu s)^2\right) \leq 2 \mu^2 s^3,
\label{ex-inverse}
\end{equation}
where the inequality completes the square in $R$ and uses $s \geq 0$. 
By Lemmas \ref{ex-origin-lemma} and \ref{ex-tail-lemma}, $A$ and $s$ tend to 0 at both ends, so $J(0)=J(\pi)=0$. 
Integrating \eqref{ex-inverse} from 0, 
which uses no sign information on $J$, 
gives $J(\sigma) \leq 2 \mu^2 \int_0^\sigma \sin ^3 t d t$.

For $0<\sigma<\pi/2$, with $c=\cos\sigma$, $I(\sigma)=\int_0^\sigma\sin^3=\frac13(1-c)^2(c+2)$, 
and $s^4-4c^2I=\frac13(1-c)^3(4c^2+9c+3)>0$; 
since $J\le2\mu^2I\le4I$ for $\mu\le\sqrt2$, $A^2=s^2+J<s^2+s^4/c^2=s^2/c^2$.
That is, $p:=A / s<\sec \sigma$ on this interval.
Next, put $T=R-2 \mu \sin \sigma$. By \eqref{ex-inverse}, $J_\sigma=-2 R s T$, 
so $J$ increases where $T<0$ and decreases where $T>0$; 
we show that $T$ changes sign at most once. Since $R_\sigma=R / A$, 
we have $T_\sigma=R / A-2 \mu \cos \sigma$, and at a zero of $T$, 
where $R=2 \mu s$, this becomes
\[
T_\sigma=2 \mu\left(\frac{1}{p}-\cos \sigma\right)>0 .
\]
Indeed, for $\sigma<\pi / 2$ positivity is the bound $p<\sec \sigma$, 
and for $\sigma \geq \pi / 2$ it is immediate because $\cos \sigma \leq 0<1 / p$. 
Thus, $T$ crosses zero only from negative to positive and transversally. 
Hence, $T$ has at most one zero: after a zero, $T$ is positive, 
so the next zero would be reached from above, forcing $T_\sigma \leq 0$ there.

We next show that $a>1 /(2 \mu)$. Write $\alpha=1 / a$. 
Inverting the origin expansion \eqref{ex-origin}, $\sigma=a r-c_3 r^3+O\left(r^5\right)$, 
gives $R(\sigma)=\alpha \sigma+c_3 \alpha^4 \sigma^3+O\left(\sigma^5\right)$, and hence
\[
T=(\alpha-2 \mu) \sigma+\left(c_3 \alpha^4+\frac{\mu}{3}\right) \sigma^3+O\left(\sigma^5\right) .
\]
Thus, when $\alpha \neq 2 \mu$, the sign of $T$ near $\sigma=0$ is that of $\alpha-2 \mu$. 
In the borderline case $\alpha=2 \mu$ we have $2 \mu a^2=a$, so $c_3=a^3 / 12$ and $c_3 \alpha^4=\mu / 6$; 
hence $T=\mu \sigma^3 / 2+O\left(\sigma^5\right)>0$.
Consequently, if $\alpha \geq 2 \mu$, then $T>0$ near $\sigma=0$, and 
since $T$ can cross zero only upward, $T>0$ on all of $(0, \pi)$. By \eqref{ex-inverse}, $J_\sigma=-2 R s T<0$ there, so 
$J$ decreases strictly from $J(0)=0$ and $J(\pi)<0$, contradicting $J(\pi)=0$. 
Therefore $\alpha<2 \mu$, that is, $a>1 /(2 \mu)$.

Since $\alpha<2 \mu$, we have $T<0$ near $\sigma=0$. 
As $\sigma \uparrow \pi$, $R(\sigma) \rightarrow \infty$ while $2 \mu \sin \sigma \leq 2 \mu$, so $T>0$ 
near $\sigma=\pi$. Hence $T$ has a zero by continuity, and 
since it has at most one, it has exactly one, at a level $\tau\in(0,\pi)$ that we call the stress transition;
moreover $T<0$ on $(0, \tau)$ and 
$T>0$ on $(\tau, \pi)$. By \eqref{ex-inverse}, $J_\sigma=-2 R s T$ is therefore positive on $(0, \tau)$ 
and negative on $(\tau, \pi)$, so $J$ increases strictly from $J(0)=0$ up to $\sigma=\tau$ 
and then decreases strictly to $J(\pi)=0$. 
Consequently $J>0$ on ($0, \pi$). Since $A, s>0$, this gives $A>\sin f$ and $A^2-\sin ^2 f>0$, which are 
the first inequality and the lower bound in the third inequality of \eqref{ex-stress-bounds}; 
the second inequality, $a>1 /(2 \mu)$, was proved above. 
Finally, as $\sin ^3 t \geq 0$, the bound $J(\sigma) \leq 2 \mu^2 \int_0^\sigma \sin ^3 t d t$ obtained at the start gives
\[
A^2-\sin ^2 f=J \leq 2 \mu^2 \int_0^\pi \sin ^3 t d t=\frac{8}{3} \mu^2,
\]
the upper bound in the third inequality.

Since $f^{\prime 2}+\sin ^2 f / r^2=\left(f^{\prime}+\sin f / r\right)^2-2 f^{\prime} \sin f / r$ 
and $-\int_0^{\infty} f^{\prime} \sin f d r=[\cos f]_0^{\infty}=2$, the 
exchange energy satisfies the exact identity (the Bogomol'nyi completion of the square, see \cite{BP75, Bog76})
\[
E_e=2+\frac{1}{2} \int_0^{\infty}\left(f^{\prime}+\frac{\sin f}{r}\right)^2 r d r .
\]
In the variable $\sigma$ we have $f^{\prime}+\sin f / r=-(A-s) / r$ and $d r=(R / A) d \sigma$, 
so with $p=A / s$,
\[
E_e-2=\frac{1}{2} \int_0^\pi \frac{(A-s)^2}{A} d \sigma=\frac{1}{2} \int_0^\pi \frac{(p-1)^2}{p} s d \sigma .
\]
By the first inequality in \eqref{ex-stress-bounds}, $A>s$, that is, $p>1$. 
Hence $(A-s)^2 / A \leq A-s$, and $(A-s)^2 \leq(A-s)(A+s)=J$ gives $A-s \leq \sqrt{J}$. 
Together with $J \leq \frac{8}{3} \mu^2$ from \eqref{ex-stress-bounds}, this yields
\[
E_e-2 \leq \frac{1}{2} \int_0^\pi \sqrt{J} d \sigma \leq \frac{\pi}{2} \sqrt{\frac{8}{3}} \mu=\pi \sqrt{2 / 3} \mu,
\]
which is the first bound in \eqref{ex-energy-bounds}. 
The inequality $A>s$ also controls $R$: 
since $(\log R)_\sigma=1 / A<1 / \sin \sigma=(\log \tan (\sigma / 2))_\sigma$, the ratio $R(\sigma) / \tan (\sigma / 2)$ 
is strictly decreasing, and by the origin expansion $R(\sigma)=\alpha \sigma+O\left(\sigma^3\right)$ 
it tends to $2 \alpha$ as $\sigma \downarrow 0$. Therefore
\begin{equation}
 R(\sigma)<2\alpha\tan(\sigma/2).                \label{ex-BP-upper}
\end{equation}
The right side is $R$ for the Belavin--Polyakov profile 
$f_B(r)=2 \arctan (2 \alpha / r)$, which has the same origin slope $a$; 
indeed, $\pi-f_B(r)=2 \arctan (r /(2 \alpha))$. 
Equivalently, \eqref{ex-BP-upper} states that $f(r)<f_B(r)$ for all $r>0$.

Next, we derive an exact identity relating $D$ and $Z$. 
Recall that
\[
J^{\prime}=2 A \sin f(2 \mu \sin f-r)=-4 \mu r f^{\prime} \sin ^2 f+2 r^2 f^{\prime} \sin f .
\]
We integrate over $(0, \infty)$. The left side gives $J(\infty)-J(0)=0$. 
On the right, the first term gives $-2 \mu D$, and since $f^{\prime} \sin f=(1-\cos f)^{\prime}$, 
integration by parts turns the second into
\[
2\left[r^2(1-\cos f)\right]_0^{\infty}-4 \int_0^{\infty} r(1-\cos f) d r=-4 Z.
\]
All boundary terms vanish by Lemmas \ref{ex-origin-lemma} and \ref{ex-tail-lemma}. Hence
\[
\mu D+2 Z=0 .
\]
This is the virial identity for the dilations $f(r) \mapsto f(r / \lambda)$, 
under which $E_e, D$ and $Z$ scale like 1, $\lambda$ and $\lambda^2$ \cite{Der64}. 
In particular $D<0$ and $2 Z=\mu|D|$. Since $f^{\prime}<0$, the Cauchy--Schwarz inequality gives
\[
|D|=2 \int_0^{\infty} r\left|f^{\prime}\right| \sin ^2 f d r \leq 2\left(\int_0^{\infty} r f^{\prime 2} d r\right)^{1 / 2}\left(\int_0^{\infty} r \sin ^4 f d r\right)^{1 / 2} \leq 4 \sqrt{Z E_e},
\]
because $\int_0^{\infty} r f^{\prime 2} d r \leq 2 E_e$, 
and $\sin ^4 f \leq \sin ^2 f=(1-\cos f)(1+\cos f) \leq 2(1-\cos f)$ gives $\int_0^{\infty} r \sin ^4 f d r \leq 2 Z$. 
Therefore $2 Z=\mu|D| \leq 4 \mu \sqrt{Z E_e}$; 
dividing by $2 \sqrt{Z}$ and squaring gives $Z \leq 4 \mu^2 E_e$. 
Finally, the first bound in \eqref{ex-energy-bounds} and $\mu \leq \sqrt{2}$ give $E_e<2+\pi\sqrt{2/3}\sqrt2<6$, so $Z<48$. 
This proves \eqref{ex-energy-bounds}.
\end{proof}

\subsection{A priori bounds and compactness}

This section shows that the shooting roots with $\mu$ in a compact subinterval of $(0,\sqrt{2}]$ form a compact set (Proposition~\ref{ex-proper}). 
Lemma~\ref{ex-stress} gives the lower bound $a>1/(2\mu)$ on the origin slope; 
Lemma~\ref{ex-steep} provides an upper bound, which rules out concentration at the origin; 
and Lemma~\ref{ex-uniform-tail} gives a uniform exponential bound at infinity. 
The proof of Lemma~\ref{ex-steep} estimates the stress $J$ at the stress transition in two ways. 
If $a\to\infty$, the transition moves to the center. 
There, a Riccati comparison for $-rf'/\sin f$ gives $J=o(\varepsilon^2)$, where $\varepsilon$ is the value of $f$ at the transition, 
whereas integrating \eqref{ex-inverse} from the center shows that $J$ is at least of order $\mu^2\varepsilon^2$. 
Since $\mu$ is bounded away from 0, the two estimates are incompatible.

\begin{lemma}[Uniform bound on the origin slope]\label{ex-steep}
For each $\mu_->0$ there is $a_*(\mu_-)<\infty$ such that
$-f'(0)\le a_*(\mu_-)$ for every $\mu\in[\mu_-,\sqrt{2}]$ and
every $f\in\Pmu$.
\end{lemma}
\begin{proof}
Suppose otherwise. Then there are $\mu_n \in[\mu_{-}, \sqrt{2}]$ and $f_n \in \mathcal{P}_{\mu_n}$ 
with $a_n=-f_n^{\prime}(0) \rightarrow \infty$, so $\alpha_n=1 / a_n \rightarrow 0$. 
Let $\tau_n$ be the stress transition of $f_n$ (see the proof of Lemma~\ref{ex-stress}), 
and let $r_n^*=R_n(\tau_n)=2 \mu_n \sin \tau_n$ be the corresponding radius. 
By \eqref{ex-BP-upper} at $\sigma=\tau_n$, $2 \mu_n \sin \tau_n<2 \alpha_n \tan (\tau_n / 2)$, 
that is, $\alpha_n>\mu_n(1+\cos \tau_n)$. 
Since $\mu_n \geq \mu_{-}>0$, this forces $\varepsilon_n=\pi-\tau_n \rightarrow 0$, and hence
\[
 r_n^*=2\mu_n\sin\varepsilon_n\to0,\qquad
 f_n(r_n^*)=\varepsilon_n\to0.
\]

We first show that $p_n(r_n^*) \rightarrow 1$, 
where $p=-r f^{\prime} / \sin f=A / s$; recall that $p>1$ by Lemma \ref{ex-stress}. 
Using the formula for $A^{\prime}$ from the proof of Lemma \ref{ex-stress}, 
direct differentiation gives
\begin{equation}
 rp'=(p^2-1)\cos f+2\mu r\sin f-r^2.            \label{ex-p}
\end{equation}
Suppose that $p_n(r_n^*) \geq 1+\delta$ for some fixed $\delta>0$ 
and infinitely many $n$. Fix $T_\delta>\log ((2+\delta) / \delta)$ and set $t=\log (r / r_n^*)$, 
so that $\partial_t=r \partial_r$. 
For $0 \leq t \leq T_\delta$ we have $r \leq r_n^* e^{T_\delta} \rightarrow 0$, 
and $f_n \leq f_n(r_n^*)=\varepsilon_n \rightarrow 0$ because $f_n$ is decreasing. 
Hence, as long as $p_n \geq 1+\delta$, so that $p_n^2-1 \geq 2 \delta$, 
dropping the nonnegative DMI term in \eqref{ex-p} gives, for all large $n$,
\[
\partial_t p_n \geq\left(p_n^2-1\right) \cos \varepsilon_n-r^2 \geq \frac{1}{2}\left(p_n^2-1\right) .
\]
In particular, $p_n$ is increasing, so it never falls below $1+\delta$, 
and this inequality holds on $[0, T_\delta]$ as long as $p_n$ exists. 
The solution of $\dot{y}=\frac{1}{2}(y^2-1)$ with $y(0)=1+\delta$ 
satisfies $(y-1) /(y+1)=\delta e^t /(2+\delta)$ 
and blows up at $t=\log ((2+\delta) / \delta)<T_\delta$. 
By comparison, $p_n \geq y$, so $p_n$ blows up at a finite positive radius. 
But $p_n=-r f_n^{\prime} / \sin f_n$ is finite at every positive radius, since $f_n$ is smooth and $0<f_n<\pi$. 
This contradiction shows that $p_n(r_n^*)<1+\delta$ for every $\delta>0$ and all large $n$, 
that is, $p_n(r_n^*) \rightarrow 1$. 
Finally, $J=s^2(p^2-1)$ and $s=\sin \tau_n=\sin \varepsilon_n$ at $\sigma=\tau_n$, so
\begin{equation}
 J_n(\tau_n)=\sin^2\varepsilon_n\bigl(p_n(r_n^*)^2-1\bigr)
                                      =o(\varepsilon_n^2). \label{ex-small-J}
\end{equation}

For every $f\in\mathcal{P}_\mu$ and $0<\sigma<\pi$,
\[
 R\le2\alpha\tan(\sigma/2),\quad
 J(\sigma)\le4\mu\alpha(1-\cos\sigma)^2,
 \quad
 0\le\frac1{\sin\sigma}-\frac1A
                         \le\frac{J}{2\sin^3\sigma}.
\]
The first bound is \eqref{ex-BP-upper}. 
By \eqref{ex-inverse}, $J_\sigma=4\mu Rs^2-2R^2s\le4\mu Rs^2$; 
inserting the first bound and integrating from $J(0)=0$ gives the second, 
because $\int_0^\sigma\tan(t/2)\sin^2t\,dt=\int_0^\sigma(1-\cos t)\sin t\,dt=(1-\cos\sigma)^2/2$. 
The third follows from $1/s-1/A=J/(As(A+s))$ and $A>s$. Since
$\int_0^{\pi/2}(1-\cos\sigma)^2/\sin^3\sigma\,d\sigma=1/2$ (put $u=\cos\sigma$),
the last two bounds give $\int_0^{\pi/2}(1/\sin\sigma-1/A)\,d\sigma\le\mu\alpha$.
As in the proof of \eqref{ex-BP-upper}, 
$\bigl(\log(R/\tan(\sigma/2))\bigr)_\sigma=1/A-1/\sin\sigma$ 
and $R/\tan(\sigma/2)\to2\alpha$ as $\sigma\downarrow0$. 
Therefore
\[
2\alpha e^{-\mu\alpha}\tan(\sigma/2)\le R(\sigma)\le2\alpha\tan(\sigma/2),\qquad 0<\sigma\le\pi/2.
\]
Integrating $J_\sigma=4\mu Rs^2-2R^2s$ over $(0,\pi/2)$ 
with the lower bound in the first term and the upper bound in the second, 
and using
$\int_0^{\pi/2}\tan(\sigma/2)\sin^2\sigma\,d\sigma=1/2$ and
$\int_0^{\pi/2}\tan^2(\sigma/2)\sin\sigma\,d\sigma=2\log2-1$, we obtain
\[
J(\pi/2)\ge4\mu\alpha e^{-\mu\alpha}-8(2\log2-1)\alpha^2.
\]
For $f_n$, since $\mu_n\in[\mu_-,\sqrt2]$ and $\alpha_n\to0$, 
we have $e^{-\mu_n\alpha_n}\to1$
and $\alpha_n/\mu_n\le\alpha_n/\mu_-\to0$, 
so $J_n(\pi/2)\ge4\mu_n\alpha_n(1-o(1))$.
By the beginning of the proof, $\alpha_n\ge\mu_n(1+\cos\tau_n)=2\mu_n\sin^2(\varepsilon_n/2)$
with $\varepsilon_n\to0$; hence $J_n(\pi/2)\ge2\mu_n^2\varepsilon_n^2(1-o(1))$.
Since $\tau_n\to\pi$, 
we have $\tau_n>\pi/2$ for large $n$, and $J_n$ increases on
$(0,\tau_n]$ by the proof of Lemma~\ref{ex-stress}. 
Thus
\[
J_n(\tau_n)\ge J_n(\pi/2)\ge2\mu_-^2\varepsilon_n^2(1-o(1)),
\]
which contradicts \eqref{ex-small-J}.
\end{proof}

Lemma~\ref{ex-tail-lemma} gives exponential decay for each profile separately, 
but the radius $R$ and the constants there depend on $f$. 
The next lemma gives a bound on the fixed region $r\ge200$ that is uniform in $\mu$ and $f$, 
as needed for the limits in Proposition~\ref{ex-proper}. 
The key input is the bound $Z<48$ from Lemma~\ref{ex-stress}, which forces $f\le1/8$ for $r\ge200$; 
there, the nonlinear term is small, 
and the comparison argument of Lemma~\ref{ex-tail-lemma} applies with explicit constants.

\begin{lemma}[Uniform exponential decay]\label{ex-uniform-tail}
For every $0<\mu \leq \sqrt{2}$, $f\in\Pmu$, and $r\ge200$,
\begin{equation}
 f(r)\le\tfrac18e^{-3(r-200)/4},\quad
 |f'(r)|\le\tfrac12e^{-3(r-200)/4},\quad
 |f''(r)|\le\tfrac13e^{-3(r-200)/4}.              \label{ex-uniform-decay}
\end{equation}
\end{lemma}
\begin{proof}
Since $f$ is decreasing, so is $1-\cos f$, 
and hence $Z\ge\int_0^r(1-\cos f(t))\,t\,dt\ge\frac{r^2}{2}\bigl(1-\cos f(r)\bigr)$. 
With $Z<48$ from \eqref{ex-energy-bounds}, at $r=200$, this gives
\[
1-\cos f(200)<\frac{96}{40000}<\frac{767}{98304}\le 1-\cos\frac18,
\]
where the last inequality is $1-\cos x\ge x^2/2-x^4/24$. 
Thus $f(200)<1/8$, and $0<f\le1/8$ on $[200,\infty)$.

On this tail, using only $r\ge1$ and $\mu \leq \sqrt{2}$, 
the three terms of $N_\mu(r,f)$ in \eqref{ex-source} are nonnegative and bounded by $\frac23f^3$, $3f^2$, 
and $\frac16f^3$, respectively; 
for the first, $f-\sin f\cos f=f-\frac12\sin2f\le\frac23f^3$, because $\sin y\ge y-y^3/6$. 
Hence
\[
0\le\frac{N_\mu(r,f)}{f}\le\frac56f^2+3f\le\frac{149}{384}<\frac{7}{16}=1-\Big(\frac34\Big)^2.
\]
Since $B_\eta=B_1-(1-\eta^2)$ and $B_1f=N_\mu(r,f)$, 
it follows that $B_{3/4}f<0$ for $r\ge200$. 
The comparison argument in Lemma~\ref{ex-tail-lemma}, 
with $R=200$ and $\eta=3/4$, gives
\[
f(r)\le f(200)\,\frac{K_1(3r/4)}{K_1(150)}\le\frac18e^{-3(r-200)/4}.
\]
For the last inequality, $K_1(x)=\int_0^\infty e^{-x\cosh t}\cosh t\,dt$ 
and $\cosh t\ge1$ give $K_1(x)\le e^{-(x-x_0)}K_1(x_0)$ for $x\ge x_0$; 
take $x=3r/4$ and $x_0=150$.

For the derivatives, recall from the proof of Lemma~\ref{ex-tail-lemma} that $(rf')'=rH$ and $rf'\to0$ at infinity, where
\[
H=\sin f\Big(\frac{\cos f}{r^2}-\frac{2\mu\sin f}{r}+1\Big).
\]
For $r\ge200$ the factor in parentheses lies between $1-1/400$ and $1+10^{-4}$, so $0<H\le2f$. 
Integrating from infinity with $\eta=3/4$ gives
\[
-rf'(r)=\int_r^\infty tH(t)\,dt\le\frac14\int_r^\infty te^{-\eta(t-200)}\,dt
=\frac14e^{-\eta(r-200)}\Big(\frac{r}{\eta}+\frac{1}{\eta^2}\Big).
\]
Dividing by $r\ge200$ gives $|f'(r)|\le\frac14\bigl(\frac43+\frac{16}{1800}\bigr)e^{-3(r-200)/4}<\frac12e^{-3(r-200)/4}$.
Finally, $f''=H-f'/r$, so
\[
|f''|\le2f+\frac{|f'|}{200}\le\Big(\frac14+\frac1{400}\Big)e^{-3(r-200)/4}
<\frac13e^{-3(r-200)/4}.
\]
\end{proof}

Lemmas~\ref{ex-stress}--\ref{ex-uniform-tail} combine into the following compactness result, 
which states that the projection $(\mu,a)\mapsto\mu$ of the set of roots is proper over $(0,\sqrt{2}]$. 
In the proof of Theorem~\ref{ex-main}, 
it is shown that the set of $\mu$ with $\mathcal P_\mu\neq\emptyset$ is closed in $(0,\sqrt{2}]$. 
In passing to the limit, the strict inequalities in \eqref{ex-class} may become non-strict; 
uniqueness for the initial value problem and the stress inequality of Lemma~\ref{ex-stress} restore them.

\begin{proposition}[Compactness of the shooting roots]\label{ex-proper}
Let $f_{\mu,a}$ be the regular origin shot and define
\[
 \mathcal R=\{(\mu,a):0<\mu \leq \sqrt{2},\ a>0,
                                      \ f_{\mu,a}\in\Pmu\}.
\]
For every $[\mu_-,\mu_+]\Subset(0,\sqrt2]$,
$\mathcal R\cap([\mu_-,\mu_+]\times(0,\infty))$ is compact.
\end{proposition}
\begin{proof}
By Lemmas~\ref{ex-stress} and~\ref{ex-steep}, 
every root $(\mu,a)$ with $\mu\in[\mu_-,\mu_+]$ lies in the rectangle $[\mu_-,\mu_+]\times[1/(2\mu_+),a_*(\mu_-)]$. 
It therefore suffices to show that the set is closed. Let $(\mu_n,a_n)$ be such roots converging to $(\bar\mu,\bar a)$. 
Then $\bar\mu\in[\mu_-,\mu_+]$ and $\bar a\ge1/(2\mu_+)>0$. 
Put $f_n=f_{\mu_n,a_n}$ and $\bar f=f_{\bar\mu,\bar a}$. 
By Lemma~\ref{ex-origin-lemma} near the origin and 
by continuous dependence on data for $r>0$, $f_n\to\bar f$ in $C^2$ on compact subsets of $[0,\infty)$. 
We show that $\bar f\in\mathcal P_{\bar\mu}$.

By construction $\bar f(0)=\pi$, and $\bar f'\le0$ as a limit of decreasing functions. 
The bounds of Lemma~\ref{ex-uniform-tail} hold for every $f_n$ with the same constants, hence also for $\bar f$; 
in particular $\bar f(r)\to0$ as $r\to\infty$. 
Next, $0\le\bar f\le\pi$. If $\bar f(r_0)\in\{0,\pi\}$ for some $r_0>0$, 
then $r_0$ is an interior extremum, so $\bar f'(r_0)=0$. 
Since the constants $0$ and $\pi$ solve the same equation, 
uniqueness for the initial value problem at $r_0$ would make $\bar f$ constant, contradicting $\bar f'(0)=-\bar a<0$. 
Thus $0<\bar f<\pi$. 
The strict stress inequality $-rf_n'>\sin f_n$ of Lemma~\ref{ex-stress} passes to the limit only as a non-strict inequality, 
but together with $0<\bar f<\pi$ it gives
\[
-r\bar f'\ge\sin\bar f>0,
\]
so $\bar f'<0$. Finally, the energy densities of $f_n$ are nonnegative and converge pointwise on $(0,\infty)$, 
and by \eqref{ex-energy-bounds}, their integrals satisfy $E_e(f_n)+Z(f_n)<54$. 
Fatou's lemma yields $E_e(\bar f)+Z(\bar f)\le54$. 
Hence $\bar f\in\mathcal P_{\bar\mu}$, that is, $(\bar\mu,\bar a)\in\mathcal R$, and the set is closed.
\end{proof}

\subsection{Nondegeneracy and local existence continuation}

This section shows that every admissible profile is nondegenerate: 
the linearization of \eqref{ex-ode} has no nonzero finite-energy solution (Lemma~\ref{ex-scalar}). 
Proposition~\ref{ex-open} uses this to continue roots to nearby parameters. 

Fix $f\in\mathcal P_\mu$ and put $s=\sin f$, $c=\cos f$, $A=-rf'$ and
\[
 Y=s+A>0,\quad
 B=2+\frac{A^2-s^2}{r^2}-\frac{2\mu}{r}s(1+c),
\]
\[
 L_f=-\partial_r^2-r^{-1}\partial_r+V_f,
 \quad V_f=\frac{\cos2f}{r^2}-\frac{2\mu\sin2f}{r}+\cos f.
\]
Here $L_f$ is the Jacobi operator, the linearization of \eqref{ex-ode} at $f$. 
The function $Y$ combines $\sin f$ with the dilation mode $-rf'$; 
for the Belavin--Polyakov profiles, these coincide and are annihilated by the exchange part of $L_f$.
In the full problem, $Y$ is no longer a zero mode but a strict supersolution: Lemma~\ref{ex-scalar} shows that $L_fY=sB>0$. 
Since $Y>0$, this makes the quadratic form of $L_f$ positive on finite-energy functions, so $L_f$ has a trivial kernel.

\begin{lemma}[Scalar supersolution and nondegeneracy]\label{ex-scalar}
For every admissible profile,
\begin{equation}
 B>1/5,\qquad L_fY=sB>s/5.                      \label{ex-barrier}
\end{equation}
If
\[
 \Hr=\left\{v:\int_0^\infty
                   (v'^2+v^2/r^2+v^2)r\,dr<\infty\right\},
\]
then $q_f[v]:=\int_0^\infty(v'^2+V_fv^2)\,r\,dr>0$ for every $v\in\Hr\setminus\{0\}$; in particular, $\ker_{\Hr}L_f=\{0\}$.
\end{lemma}
\begin{proof}
Differentiating \eqref{ex-ode} gives
\[
 L_fs=s\left(\frac{A^2-s^2}{r^2}-\frac{2\mu sc}{r}\right),
 \qquad L_f(rf')=2s(\mu s/r-1);
\]
for the second identity, 
note that $f(r/\lambda)$ solves \eqref{ex-ode} with $2\mu$ and $\sin f$ replaced by $2\mu/\lambda$ and $\lambda^{-2}\sin f$, 
and differentiate at $\lambda=1$. 
Hence $L_fY=sB$. 
With dots denoting $r\partial_r$, we have, using \eqref{ex-ode} for $\dot A$,
\[
 \dot s=-Ac,\quad \dot c=As,\quad
 \dot A=-sc+2\mu rs^2-r^2s,
\]
\begin{equation}
 \dot B=4-2B-2As+\frac{2\mu}{r}(1+c)(A-s).      \label{ex-Bdot}
\end{equation}
At a point where $B=1/5$, completion of the square gives $A^2\le s^2+\frac59\mu^2s^2(1+c)^2$; 
hence, as $\mu^2\le2$, $(As)^2\le s^4+\frac{10}{9}s^4(1+c)^2\le1+\frac{10}{9}\cdot\frac{1024}{729}<\frac{81}{25}$, 
since $\max_{[-1,1]}(1-c^2)^2(1+c)^2=\frac{1024}{729}$ at $c=1/3$. 
Thus $As<9/5$. Since $A>s$ by Lemma~\ref{ex-stress}, the last term in \eqref{ex-Bdot} is nonnegative, so $\dot B\ge\frac{18}5-2As>0$.
The origin expansion gives $B=2-a^2r^2/2+O(r^4)$. If
$B\le1/5$ somewhere, there is a smallest $r_1>0$ with $B(r_1)=1/5$; 
then $B>1/5$ on $(0,r_1)$, so $\dot B(r_1)\le0$, a contradiction. 
This proves \eqref{ex-barrier}.

For $v\in C_c^\infty(0,\infty)$, writing $v=Yw$ and integrating by parts gives
\begin{equation}
 \int(v'^2+V_fv^2)r\,dr
 =\int r\left(v'-\frac{Y'}Yv\right)^2dr
       +\int\frac{L_fY}{Y}v^2r\,dr.              \label{ex-groundstate}
\end{equation}
Lemmas~\ref{ex-origin-lemma} and~\ref{ex-tail-lemma} imply $Y=2ar+O(r^3)$ at zero and $Y=c_frK_1(r)(1+o(1))$ at infinity, 
with differentiated estimates. 
Thus $|Y'/Y|\le C_f(1+1/r)$. Also, $L_fY/Y=sB/Y$ is positive and bounded: it tends to $1$ at zero and is asymptotic to $2/r$ at infinity. 
Since also $|V_f|\le C_f(1+r^{-2})$, all three quadratic forms in \eqref{ex-groundstate} are continuous on $\mathcal H_r$. 
Moreover, $C_c^\infty(0,\infty)$ is dense in $\mathcal H_r$: 
if $\chi$ cuts off near $0$ and $\infty$ with $|\chi'|\le C/r$, 
then the $\mathcal H_r$-norm of $v-\chi v$ is bounded by the tails of $\int(v'^2+v^2/r^2+v^2)\,r\,dr$, 
which tend to zero, and $\chi v$ is approximated by mollification on a compact annulus. 
Hence \eqref{ex-groundstate} holds for every $v\in\mathcal H_r$, and by \eqref{ex-barrier} its second term is at least
\[
\frac15\int\frac{\sin f}{\sin f-rf'}\,v^2\,r\,dr>0\qquad\text{if }v\ne0;
\]
this proves $q_f[v]>0$.
If $v\in\mathcal H_r$ is a weak solution of $L_fv=0$, 
then by density we may test the equation with $v$ itself; 
the left side of \eqref{ex-groundstate}, which is $q_f[v]$, vanishes, so $v=0$.
\end{proof}

Proposition~\ref{ex-open} continues a root by matching the origin shot with a decaying tail at a fixed radius $R$; 
the next lemma describes these tails. 
By \eqref{ex-source}, equation \eqref{ex-ode} reads $B_1u=N_\mu(r,u)$ on $[R,\infty)$. 
Near the zero solution, its decaying solutions form a one-parameter family indexed by $p=u(R)$, 
and in the $(u(R),u'(R))$-plane they form the graph of a smooth function $\Phi_{\mu,R}$, 
the stable manifold of the zero solution at infinity. 
The proof of Lemma~\ref{ex-tail-lemma} showed that a given admissible profile satisfies the integral equation 
$f=f(R)K_1/K_1(R)+T_RN_\mu(f)$, with $R$ depending on $f$. 
Here we solve the same equation for every small boundary value $p$ at a fixed radius, 
using a contraction argument in a weighted norm with the decay rate $3/4$ from Lemma~\ref{ex-uniform-tail}; 
this gives existence, uniqueness, and smooth dependence on $(\mu,p)$, as needed for the implicit function theorem.
By Lemma~\ref{ex-uniform-tail}, every admissible profile lies on this graph at $R=204$, 
so near a root, a shot from the origin is admissible exactly when $f'(R)=\Phi_{\mu,R}(f(R))$; 
this reduces the matching to one scalar equation.

\begin{lemma}[Small decaying solutions]\label{ex-tailgraph}
Fix $R\ge1$. For $0<\mu \leq 3/2$ and $|p|\le1/128$ there is a unique
solution of \eqref{ex-ode} on $[R,\infty)$ with
\[
 u(R)=p,\qquad
 \|u\|_{\eta,R}:=\sup_{r\ge R}e^{\eta(r-R)}|u(r)|\le1/64,
 \qquad\eta=3/4.
\]
Its boundary derivative is a smooth function $u'(R)=\Phi_{\mu,R}(p)$.
For $p>0$ this solution is positive, strictly decreasing and has finite
exchange and Zeeman energy on the tail. For $0<\mu\le\sqrt2$, every admissible profile is in
this chart at $R=204$: for $f\in\mathcal P_\mu$, 
the restriction of $f$ to $[204,\infty)$ is the solution with $p=f(204)$, so $f'(204)=\Phi_{\mu,204}(f(204))$.
\end{lemma}
\begin{proof}
By \eqref{ex-source}, equation \eqref{ex-ode} reads $B_1u=N_\mu(r,u)$. 
With the Green operator $T_R$ of \eqref{ex-green}, consider
\begin{equation}\label{ex-tail-fp}
u=p\,\frac{K_1(r)}{K_1(R)}+T_RN_\mu(r,u).
\end{equation}
Every solution of \eqref{ex-tail-fp} with $\|u\|_{\eta,R}<\infty$ solves \eqref{ex-ode} with $u(R)=p$. 
Conversely, if $u$ solves \eqref{ex-ode} on $[R,\infty)$ with $u(R)=p$ and $\|u\|_{\eta,R}<\infty$, 
then the difference of the two sides of \eqref{ex-tail-fp} solves $B_1v=0$, vanishes at $R$ and tends to zero, 
so it vanishes, as in the proof of Lemma~\ref{ex-tail-lemma}. 
It therefore suffices to solve \eqref{ex-tail-fp} in the ball $\|u\|_{\eta,R}\le d:=1/64$.

Since
\[
B_1e^{-\eta(r-R)}=\big(1-\eta^2+\eta/r+r^{-2}\big)e^{-\eta(r-R)}\ge\tfrac7{16}\,e^{-\eta(r-R)},
\]
the function $\frac{16}{7}\|g\|_{\eta,R}\,e^{-\eta(r-R)}$ is a supersolution of $B_1w=|g|$, 
and positivity of $G_R$ and comparison give $\|T_Rg\|_{\eta,R}\le\frac{16}{7}\|g\|_{\eta,R}$.
For $|u|\le d$ and $r\ge R\ge1$,
\[
|\partial_uN_\mu(r,u)|=\Big|\frac{2\sin^2u}{r^2}+\frac{2\mu\sin2u}{r}+1-\cos u\Big|
\le4\mu d+\tfrac52d^2\le6d+\tfrac52d^2=\tfrac{773}{8192}.
\]
Since $|u(r)|\le\|u\|_{\eta,R}$, the right side of \eqref{ex-tail-fp} is a contraction on the ball, 
with constant $\frac{16}{7}\cdot\frac{773}{8192}=\frac{773}{3584}<\frac12$. 
It maps the ball into itself: $N_\mu(r,0)=0$, $|p|\le d/2$ 
and $K_1(r)/K_1(R)\le e^{-(r-R)}$ (as in the proof of Lemma~\ref{ex-uniform-tail}), 
so the image has norm at most $d/2+\frac{773}{3584}d<d$. 
The contraction theorem gives the unique solution, 
which depends smoothly on $(\mu,p)$ because the map does. 
Since $H_R'(R)=1/R$ by the Wronskian identity, 
differentiating \eqref{ex-tail-fp} at $r=R$ gives the convergent expression
\[
\Phi_{\mu,R}(p)=p\,\frac{K_1'(R)}{K_1(R)}+\frac{1}{RK_1(R)}\int_R^\infty tK_1(t)N_\mu(t,u(t))\,dt .
\]
Differentiating \eqref{ex-tail-fp} in $r$ 
and in the parameters likewise gives exponential bounds for $u'$ and for the parameter derivatives of $u$.

Let $p>0$. Since $G_R\ge0$ and $N_\mu(r,u)\ge0$ for $u\ge0$, 
every iterate starting from $pK_1/K_1(R)$ is at least $pK_1/K_1(R)$, so $u\ge pK_1/K_1(R)>0$. 
Moreover, by \eqref{ex-ode},
\[
(ru')'=\sin u\Big(\frac{\cos u}{r}+r-2\mu\sin u\Big)>0,
\]
because $r\ge1$, $\mu\le3/2$ and $0<u\le1/64$.
Thus $ru'$ is strictly increasing and tends to zero at infinity, so it is negative; 
that is, $u'<0$. 
The exponential decay of $u$ and $u'$ gives finite exchange and Zeeman energy on the tail.

Finally, let $f\in\mathcal P_\mu$. By \eqref{ex-uniform-decay},
\[
\sup_{r\ge204}e^{3(r-204)/4}f(r)\le e^{-3}/8<1/128 .
\]
Thus $p=f(204)<1/128$ and $\|f\|_{\eta,204}<d$, 
so the restriction of $f$ to $[204,\infty)$ is the solution above with $R=204$.
\end{proof}

The next proposition is the open half of the open--closed argument for Theorem~\ref{ex-main}; 
Proposition~\ref{ex-proper} yields the closed half, and Lemma~\ref{ex-seed} guarantees nonemptiness. 
Its proof is a shooting argument. 
Near a root, a shot $f_{\mu,a}$ from the origin is admissible exactly when its data at $R=204$ lie on the tail graph of Lemma~\ref{ex-tailgraph}; 
that is, when the matching function $S(\mu,a)=f'_{\mu,a}(R)-\Phi_{\mu,R}(f_{\mu,a}(R))$ vanishes. 
The implicit function theorem then continues the root to nearby $\mu$, provided $\partial_aS\neq0$. 
This transversality holds because a tangency would turn the $a$-derivative of the shot into a nonzero finite-energy Jacobi field, 
which Lemma~\ref{ex-scalar} excludes.

\begin{proposition}[Openness]\label{ex-open}
The set
$\mathcal U=\{\mu\in(0,\sqrt2]:\Pmu\ne\varnothing\}$ is open in
$(0,\sqrt{2}]$.
\end{proposition}
\begin{proof}
Let $(\mu_0,a_0)\in\mathcal R$, put $f_0=f_{\mu_0,a_0}$ and $R=204$. 
By Lemma~\ref{ex-tailgraph}, $f_0'(R)=\Phi_{\mu_0,R}(f_0(R))$, 
and by  \eqref{ex-uniform-decay}, $0<f_0(R)\le e^{-3}/8<1/128$. 
Near $(\mu_0,a_0)$, define
\[
S(\mu,a)=f'_{\mu,a}(R)-\Phi_{\mu,R}\big(f_{\mu,a}(R)\big).
\]
By continuity, $0<f_{\mu,a}(R)<1/128$ in a small neighborhood of $(\mu_0,a_0)$, so $S$ is defined there; 
it is smooth by Lemmas~\ref{ex-origin-lemma} and~\ref{ex-tailgraph}, and $S(\mu_0,a_0)=0$. 
Since Lemma~\ref{ex-tailgraph} covers $\mu\le3/2$, this neighborhood is open in $\R^2$ also when $\mu_0=\sqrt2$.

Shrinking the neighborhood, 
we may also assume $0<f_{\mu,a}<\pi$ and $f'_{\mu,a}<0$ on $(0,R]$: on $(0,r_0]$ for a small $r_0>0$, 
this follows from the uniform origin expansion of Lemma~\ref{ex-origin-lemma}, 
and on $[r_0,R]$ from continuous dependence, 
since $f_0$ stays away from $0$ and $\pi$ and $f_0'<0$ there. 
If $S(\mu,a)=0$, then $f_{\mu,a}$ and the solution $u$ of Lemma~\ref{ex-tailgraph} with $p=f_{\mu,a}(R)$ have the same value and derivative at $R$, 
so they coincide on $[R,\infty)$ by uniqueness for the initial value problem. 
By Lemma~\ref{ex-tailgraph}, this tail is positive, strictly decreasing, decaying, and of finite energy, so $f_{\mu,a}\in\mathcal P_\mu$. 
Conversely, every root near $(\mu_0,a_0)$ satisfies $S=0$ by Lemma~\ref{ex-tailgraph}. 
Thus, near $(\mu_0,a_0)$, the roots are exactly the zeros of $S$.

It remains to show that $S_a(\mu_0,a_0)\neq0$; 
the implicit function theorem then gives a root $(\mu,a(\mu))$ for every $\mu$ near $\mu_0$, so $\mathcal U$ is open. 
Let $v=\partial_af_{\mu,a}|_{(\mu_0,a_0)}$ and $z=\partial_pu|_{p=f_0(R)}$, 
where $u$ is the tail of Lemma~\ref{ex-tailgraph} at $\mu=\mu_0$; thus $u=f_0$ on $[R,\infty)$ when $p=f_0(R)$. 
Differentiating \eqref{ex-ode} gives $L_{f_0}v=0$ on $(0,R]$ and $L_{f_0}z=0$ on $[R,\infty)$, and
\[
v=-r+O(r^3),\qquad z(R)=1,\qquad z'(R)=\Phi_p,\qquad S_a=v'(R)-\Phi_p\,v(R),
\]
where $\Phi_p=\partial_p\Phi_{\mu_0,R}(f_0(R))$. 
Here $(v(R),v'(R))$ is the velocity of the shot data at $R$ as $a$ varies, 
and $(1,\Phi_p)$ is tangent to the tail graph; $S_a=0$ means that they are parallel. 
In that case let $w=v$ on $[0,R]$ and $w=v(R)z$ on $[R,\infty)$. 
The values and derivatives of the two pieces agree at $R$, so $w$ solves $L_{f_0}w=0$ on $(0,\infty)$. 
Moreover, $w\in\mathcal H_r$ due to the origin expansion 
and the exponential bounds from the proof of Lemma~\ref{ex-tailgraph}, and $w\neq0$ because $w=-r+O(r^3)$. 
This contradicts Lemma~\ref{ex-scalar}. Hence $S_a(\mu_0,a_0)\neq0$.
\end{proof}

\subsection{Existence at \texorpdfstring{$\mu=1/2$}{mu = 1/2} and proof of Theorem~\ref{ex-main}}

The open--closed argument needs one parameter at which an admissible profile exists. 
We take it from Ibrahim--Shimizu~\cite{IS}: after an exact dilation, 
their model with $\beta=1$ coincides with ours, 
and their existence and monotonicity results give an admissible profile at $\mu=1/2$.

\begin{lemma}[One nonempty fiber]\label{ex-seed}
$\mathcal P_{1/2}\ne\varnothing$.
\end{lemma}
\begin{proof}
We use the following restricted content of Ibrahim--Shimizu \cite[Theorem~A and equation (1.13)]{IS}. 
In their radius $\rho$ and 
coupling $\mathfrak r$ (we write $\mathfrak r$ for their coupling to avoid a clash with our radius), 
at $\beta=1$ and $0<\mathfrak r\le1$ there
is a profile $F$ with endpoints $\pi,0$, $0<F<\pi$, finite radial
exchange, and $F\in L^2(\rho\,d\rho)$, satisfying
\[
 -F_{\rho\rho}-F_\rho/\rho+\sin F\cos F/\rho^2
       -2\mathfrak r\sin^2F/\rho+2\sin F=0,
 \qquad F_\rho+\sin F/\rho<0.
\]
Its Bloch helicity and curl convention are the ones used here.
Set $\mathfrak r=\sqrt2\mu$ and $f(r)=F(r/\sqrt2)$.
Multiplying this displayed equation by $1/2$ gives exactly
\eqref{ex-ode} for $f$. 
In \cite{IS} the equation holds in the weak sense; since its coefficients are smooth on $(0,\infty)$, 
a bootstrap gives $F\in C^\infty(0,\infty)$, and hence $f\in C^\infty(0,\infty)$.
Exchange is unchanged by this dilation, and
\[
 Z(f)=2\int_0^\infty(1-\cos F)\rho\,d\rho
                         \le\int_0^\infty F^2\rho\,d\rho<\infty.
\]
For $\mu=1/2$ the required coupling is $1/\sqrt2<1$, so all the
quoted hypotheses hold. Thus $f\in\mathcal P_{1/2}$.
\end{proof}

\begin{proof}[Proof of Theorem~\ref{ex-main}]
By Proposition~\ref{ex-open}, the image $\mathcal U$ of the root
projection is open in $(0,\sqrt{2}]$. It is also relatively closed: if
$\mu_n\in\mathcal U$ tends to $\mu_*\in(0,\sqrt2]$, choose any root at
each $\mu_n$. Eventually the parameters lie in a compact subinterval
of $(0,\sqrt{2}]$, and Proposition~\ref{ex-proper} gives a limiting root
at $\mu_*$. Lemma~\ref{ex-seed} makes $\mathcal U$ nonempty.
Connectedness of $(0,\sqrt{2}]$ therefore gives $\mathcal U=(0,\sqrt{2}]$.
The endpoint conclusions were proved for every admissible profile in
Lemmas~\ref{ex-origin-lemma} and~\ref{ex-tail-lemma}.
\end{proof}

\begin{remark}
After the dilation in the proof of Lemma~\ref{ex-seed}, \cite[Theorem~A]{IS} gives $\mathcal P_\mu\neq\emptyset$ for every $0<\mu\le1/\sqrt2$. 
We use it only at $\mu=1/2$; all other parameters are reached by the continuation argument.
\end{remark}

\end{document}